# Inter-Stop Energy Prediction and Causal Driver Quantification for Dual-Source Trolleybuses via a Time-Aware Tabular Deep Learning Architecture

Wentao Zeng [a,b] , Zijian Huang [c], Yiming Bie [d], Jiabin Wu [a,*], Jun Gong [e]

[a] *School of Management, Foshan University, Foshan 528000, China*
[b] *School of Mechanical and Electrical Engineering and Automation, Foshan University, Foshan 528225, China*
[c] *School of Artificial Intelligence, South China Normal University, Guangzhou 510631, China*
[d] *School of Transportation, Jilin University, Changchun 130022, China*
[e] *Department of Civil Engineering, The University of Hong Kong, Hong Kong 999077, China*

* **Corresponding author: Jiabin Wu** (E-mail address: jiabinwu@fosu.edu.cn)

**Abstract:** Dual-source trolleybuses switch between overhead catenary power and on-board battery operation, producing energy patterns shaped by static route attributes, high-frequency vehicle trajectories, and hourly meteorological conditions. Existing models encode such heterogeneous inputs inadequately, and none provides causal explanation of the factors governing energy use. This paper addresses both limitations. On the artificial intelligence side, a time-aware tabular deep learning architecture is proposed. Periodic time encoding is embedded into a parameter-efficient batch-ensemble backbone, enabling static and sequential features to be represented jointly in a single network. Bayesian optimization guided by tree-structured density estimation drives the hyperparameter search. A three-layer causal explanation pipeline then extends the predictive model: feature-level attribution quantifies marginal contributions, a linear non-Gaussian acyclic model recovers causal directions from observational data, and a meta-learner estimates net average treatment effects — moving the analysis from correlation to causal mechanism. On the engineering side, the framework targets inter-stop energy management for dual-source trolleybuses in urban transit networks. Experiments on the Zurich trolleybus dataset, enriched with meteorological records, yield a mean absolute percentage error of 6.52% and a coefficient of determination of 0.982, surpassing all ten baseline models across statistical, tree-ensemble, and deep learning categories. Ablation experiments confirm periodic time encoding as the principal source of accuracy gain. Causal analysis identifies the regenerative braking ratio and average speed as the dominant energy-saving levers. Coasting distance, by contrast, emerges as the primary contributor to excess consumption. These results provide concrete intervention thresholds for vehicle technology, driving behavior, capacity allocation, and catenary network planning.

# 1. Introduction

Transport has become the largest urban source of air pollution in most cities worldwide, and pressure for energy saving and emission reduction continues to mount [1-3]. The IEA Global EV Outlook 2025 [4] reports that in 2024 the global electric fleet displaced more than 1.3 million barrels of oil per day, and electric bus sales grew by 30% year-on-year to surpass 70,000 units, underscoring the central role of transit electrification in the energy transition. Pure battery-electric technology, however, still faces a structural payload-range-cost dilemma in heavy-duty commercial applications: highly redundant battery packs account for close to 50% of the total purchase cost, which severely constrains operational economics. Although trams, trolleybuses, and battery-electric buses have all been deployed at scale, reconciling operational efficiency and range flexibility under limited infrastructure investment remains a critical bottleneck for the sector.

The dual-source trolleybus (DTB), with its hybrid supply architecture, offers a competitive solution to this dilemma. In grid-covered sections, the DTB draws power dynamically from an overhead catenary through a pantograph, reducing reliance on large batteries; in off-grid sections, it runs on an on-board battery and is no longer bound to a fixed catenary network as conventional trolleybuses are. The key advantage of this architecture is that a DTB does not need an oversized battery for the worst-case range nor full-route catenary coverage: a localized network at key segments is sufficient to jointly optimize vehicle mass, range flexibility, and deployment cost [5]. The energy behavior of a DTB, however, inherits the complexity of both battery-driven and catenary-fed operation. Energy use is shaped not only by vehicle parameters but also by route grade, stop frequency, passenger load, catenary coverage, and weather conditions. Without a systematic understanding of how these factors drive energy use, operators can neither properly evaluate catenary layout at the planning stage nor design differentiated energy-management strategies at the operational stage. Improving the accuracy of DTB energy prediction and uncovering its underlying mechanisms is therefore of both theoretical and practical importance for fine-grained fleet management.

Building a high-accuracy and interpretable energy model for DTBs raises three core challenges. (1) Heterogeneous data fusion and spatio-temporal feature extraction are inherently difficult. DTB energy consumption arises from nonlinear couplings among vehicle attributes, network topology, driving behavior, and environmental factors. These couplings yield a heterogeneous feature space that combines two feature types. Static attributes capture section length and grade, while high-frequency dynamic series capture speed trajectories and temperature variation. The two types differ markedly in sampling frequency, units, and distribution. Physics-based models are interpretable but depend on empirically calibrated parameters and struggle to absorb random traffic disturbances, while deep sequence models handle temporal dynamics well but remain less efficient than gradient-boosted trees on sparse categorical features. A unified architecture that combines efficient tabular encoding with sequential modeling is therefore needed. (2) Deep models are sensitive to hyperparameters, and optimization efficiency is a bottleneck. In Time-Aware Tabular DL model that makes Multiple predictions (TabM), the hyperparameter space spans learning rate, hidden-layer width, ensemble size and other dimensions that are nonlinearly coupled. Grid search is computationally prohibitive [6], while random search lacks adaptivity [7]; approaching the global optimum under a limited compute budget thus becomes a key issue. (3) Predictive models are typically opaque and offer no causal explanation. Most high-accuracy energy models are black boxes [8]: they may observe that low speed correlates with high energy use, when the true driver is traffic congestion that causes both. Conflating correlation with causation severely limits the translation of predictions into operational decisions. Turning such black-box algorithms into transparent planning tools [9] is a pressing issue for intelligent bus energy management.

To address these gaps, this study starts from multi-source data fusion and feature engineering, integrates micro-level vehicle state with macro-level meteorological data, and constructs a high-quality heterogeneous dataset that contains both static and dynamic features through spatio-temporal alignment and ratio-based feature engineering. A Time-Aware TabM architecture is then built by embedding Time2Vec periodic encoding into the TabM tabular backbone, unifying tabular feature encoding and temporal modeling in a single network to address challenge (1). Tree-structured Parzen Estimator (TPE)-based Bayesian optimization with median pruning is used to tune the model, enabling efficient convergence toward the global optimum under a limited compute budget and addressing challenge (2). Shapley Additive exPlanations (SHAP) is used to quantify marginal feature contributions, Direct Linear Non-Gaussian Acyclic Model (DirectLiNGAM) recovers a causal directed acyclic graph (DAG) over the variables, and a Two-Learner (T-Learner) meta-algorithm [10] estimates the net average treatment effect (ATE) of key drivers, moving the analysis from correlation to causation and addressing challenge (3). The main contributions of this paper are as follows.

(1) Focusing on inter-stop energy modeling for hybrid-supplied DTBs, this study places dual-mode supply behavior in a unified analytical view, resolves the mismatch between second-level telemetry and hourly meteorological records through spatio-temporal alignment and ratio-based feature engineering, and builds a high-quality modeling dataset with 63 fields grouped into 8 categories.

(2) This study proposes Time-Aware TabM, a unified predictive architecture that integrates temporal awareness into tabular modeling. Existing approaches typically encode temporal and tabular features in separate branches. In contrast, the proposed architecture adopts the parameter-efficient TabM ensemble as its backbone and embeds Time2Vec encoding, which is coupled with a flattened GEMM batch-ensemble. Within a single network, static tabular features and dynamic temporal features are therefore modeled jointly.

(3) A three-layer causal explanation system — "SHAP attribution → DirectLiNGAM causal discovery →T-Learner effect estimation" — is constructed. Unlike previous work that relies on a single post-hoc explanation method, this study integrates structural causal discovery and counterfactual effect estimation into the analysis of DTB energy use, moving from feature importance ranking to causal direction identification and then to quantitative intervention effects.

(4) The framework is validated through ablation studies and benchmarked against ten baseline models. Causal-inference results then identify two core intervention anchors: the regenerative braking ratio and the average speed. Building on these anchors, actionable operational recommendations are proposed with concrete parameter thresholds. The recommendations span four dimensions. These dimensions cover vehicle technology, driving behavior, capacity allocation, and catenary-network planning.

The remainder of the paper is organized as follows. Section 2 reviews the related literature. Section 3 describes the data sources and preprocessing. Section 4 presents the methodology. Section 5 reports the experimental results and discussion. Section 6 concludes the paper and outlines future work.

# 2. Literature Review

Research on electric-bus energy prediction has advanced along two broad lines: physics- or statistics-based approaches, and data-driven approaches. Despite substantial progress, both lines leave unresolved gaps in heterogeneous feature fusion, temporal modeling, and causal interpretation—gaps that become particularly acute for the dual-source trolleybus (DTB), whose hybrid catenary–battery supply architecture introduces additional complexity absent from conventional battery-electric bus studies. Model interpretability and causal inference have recently emerged as a third frontier, yet existing work in this area has yet to be integrated with high-accuracy predictive models for hybrid-supplied vehicles. This section reviews the literature along these three dimensions, identifying the specific limitations that motivate the present work; Appendix Table A-1 provides a structured contrast of representative studies from the past five years against the proposed framework.

## 2.1 Physics- and statistics-based energy prediction

Physics-based models reconstruct energy consumption by decomposing vehicle dynamics into discrete power pathways—traction, braking, and auxiliary loads—and mapping external operating conditions to the evolving battery state of charge. This mechanistic transparency has established them as an enduring reference standard in electric-bus modeling. At the microscopic level, Fiori et al. [11] demonstrated that decoupling traction power from auxiliary power yields fine-grained resolution of energy flows across heterogeneous operating conditions, a modeling template that subsequent studies have built upon. Hjelkrem et al. [12] scaled this approach to the system level by coupling longitudinal vehicle dynamics with HVAC airflow exchange, thereby incorporating the influence of ambient temperature on auxiliary systems into a unified white-box framework and substantially improving physical consistency for urban bus fleets across diverse climates. Extending the paradigm to real-time use, Pan et al. [13] paired longitudinal dynamics with recursive least squares incorporating a forgetting factor, enabling stop-by-stop cumulative energy prediction and demonstrating that physics-based models can be made adaptive through online parameter estimation. At the planning level, Lajunen [14] conducted vehicle-dynamics simulations across multiple route profiles, producing a cost–benefit assessment of hybrid and battery-electric buses that remains a useful reference for fleet procurement. Most recently, Ferdaus et al. [15] surveyed the zero-shot potential of foundation models for clean-energy forecasting; while promising, their domain-specific adaptation to DTB telemetry has not yet been demonstrated.

Physics-based models thus offer clear mechanistic interpretability and transparent derivations. Their structural strength, however, is inseparable from a critical dependency on accurately calibrated parameters.

Under the compound complexity of urban DTB operation—where stochastic driving behavior, volatile traffic flow, and multi-source energy coupling interact simultaneously—this parameter sensitivity translates directly into accuracy degradation and limited generalizability. The inability of white-box models to absorb unmodeled disturbances without manual recalibration has been a primary impetus for the rise of data-driven methods, and it defines the modeling gap that the present framework is designed to bridge.

## 2.2 Data-driven energy prediction

The widespread deployment of vehicle telematics and smart transit infrastructure has generated large-scale, multi-source operational records, enabling data-driven methods to learn energy patterns directly from trajectories without relying on physical parameter calibration. Environmental influences were among the first to benefit from this shift. Kim et al. [16] applied regression models to auxiliary-energy prediction under extreme weather, quantifying nonlinear thermal loads and establishing ambient temperature as a primary driver of non-traction energy use. Li et al. [17] combined a stochastic speed-profile generator with a KNN–Random Forest (KNN–RF) hybrid, using ablation decomposition to show that grade–congestion interaction has been systematically underestimated in single-factor analyses—a finding with direct implications for route-level energy management. The scope of data-driven modeling then widened to driving behavior and multi-factor dynamics. Zhang et al. [18] constructed an XGBoost-based three-level driving-segment framework on large-scale Beijing taxi records, enabling energy analysis spanning macroscopic traffic flow to microscopic driving events, though the taxi setting limits transferability to transit bus operations. Zhang et al. [19] pursued a nested optimization approach on fleet data to characterize battery-electric bus energy patterns and inform charging-strategy design. Liu and Song [20] embedded energy prediction within a joint planning model for dynamic wireless-charging siting, linking prediction quality directly to infrastructure decisions. For online deployment, Chen et al. [21] showed that Long Short-Term Memory (LSTM) and Artificial Neural Network (ANN) with a data-reconstruction algorithm can capture instantaneous energy fluctuations in real time, while Pan et al. [22] introduced an incremental-learning scheme that continuously adapts model parameters to evolving operating conditions—yet both works treat hyperparameter selection as secondary and offer no mechanism for causal explanation.

Collectively, these works demonstrate that data-driven methods can attain high predictive accuracy across a range of bus types and operating environments. Two structural deficiencies, however, persist across virtually all of them. First, none addresses the principled joint encoding of static tabular attributes and high-frequency dynamic series within a single architecture—a requirement that becomes acute for DTBs, whose inter-stop feature space spans route geometry, catenary-coverage indicators, second-level speed trajectories, and hourly meteorological records simultaneously. Feeding such heterogeneous inputs into sequence models designed for uniform time series, or into tree ensembles that discard temporal order, produces architectures that are poorly matched to the data structure. Second, the dominant emphasis on prediction accuracy leaves the mechanisms of energy use unexplained: a model that learns low speed correlates with high consumption cannot determine whether congestion is the common cause, or whether the relationship reverses under regenerative braking. This confusion between statistical association and causal direction severely constrains the operational utility of predictions. These two gaps collectively define the design space for the time-aware predictive architecture and causal inference pipeline proposed in the present study.

### 2.3 Model interpretability and causal inference

Growing dissatisfaction with opaque predictive models has driven a parallel research strand focused on interpretability and causal inference. In transit applications, Zhu et al. [23] developed a causal AI framework employing a conditional correction index to overcome correlation bias in carbon-efficiency assessment, marking an early demonstration that causal reasoning can improve the policy relevance of transit analytics. Most directly related to the present study, He et al. [24] paired a Transformer backbone with double/debiased machine learning to estimate causal effects of inter-stop energy use for electric buses under sparse data conditions; while this work shares our focus on stop-level causal analysis, it employs a fundamentally different backbone architecture without time-awareness encoding and does not address the heterogeneous tabular–temporal fusion problem central to DTB modeling. Llopis-Castelló et al. [25] used a before-after design with connected-vehicle data to quantify energy and emission benefits from pavement-roughness interventions, illustrating the operational value of causal evaluation in infrastructure management. On the interpretability side, Song et al. [26] combined a Kolmogorov–Arnold network with multi-scale dynamic feature normalization to improve the transparency of deep models in power-load forecasting—a non-transport domain where the absence of driving-behavior features limits direct transferability. Zhang et al. [27] applied double/debiased difference-in-differences estimation to quantify the social impact of the London Night Tube, establishing a rigorous template for causal inference in transport policy evaluation. Jia et al. [28] advanced spatial interpretability through a multi-scale deep learning model for urban energy demand forecasting, while Eskandari et al. [29] integrated feature selection with interpretable machine learning at the national energy scale; neither framework is designed for stop-level heterogeneous data. Zhang et al. [30] constructed a spatio-temporal causal game model to disentangle heterogeneous correlations across charging stations, demonstrating that causal discovery can inform infrastructure-level load balancing—though the charging-station context is structurally distinct from the in-service bus energy prediction task considered here.

These studies collectively advance model transparency and establish causal inference as a viable methodology in transit research. Yet two structural gaps persist. First, the majority of interpretability-focused works rely on a single post-hoc attribution tool—typically SHAP values or attention weights—and do not integrate structural causal discovery or quantitative effect estimation into the analytical pipeline; attribution rankings and causal directions consequently remain conflated. Second, and more specifically, no prior study has applied a unified causal framework to hybrid-supplied DTBs, whose dual-mode supply architecture creates energy interdependencies that conventional bus energy studies do not encounter. The present study addresses both gaps by constructing a three-layer progressive explanation system: SHAP multi-layer attribution quantifies feature-level marginal contributions, DirectLiNGAM recovers the causal directed acyclic graph from observational data, and T-Learner meta-estimation converts causal directions into net average treatment effects—forming a coherent chain from correlation ranking to mechanism quantification.

The review above reveals two tightly coupled open problems in intelligent transit energy management. How to jointly encode heterogeneous tabular and temporal features within a single architecture capable of handling DTB-specific data structures, and how to introduce causal inference that is integrated with—rather than merely appended to—a high-accuracy predictive model. Existing approaches resolve at most one of these problems at a time: physics-based models offer interpretability but sacrifice flexibility; data-driven models achieve accuracy but remain opaque; and interpretability-focused works apply causal tools without the high-fidelity predictive backbone needed to produce reliable counterfactuals. The present study addresses both problems in concert. A Time-Aware TabM architecture embeds Time2Vec periodic encoding into a parameter-efficient TabM backbone, unifying tabular and temporal feature modeling under a single training objective. This predictive foundation then supports a three-layer causal explanation pipeline—DirectLiNGAM structural discovery, SHAP multi-layer attribution, and T-Learner treatment-effect estimation—that moves the analysis of DTB inter-stop energy use from statistical association to quantified causal mechanism.

## 3. Data Collection and Preprocessing

The full data collection and preprocessing pipeline is illustrated in Fig. 1. Two heterogeneous data sources are integrated. The first is the ZTBus benchmark dataset [31]. In terms of data acquisition, signals from on-board components—such as the motor, wheel encoders, and temperature sensors—are streamed to the Vehicle Control Unit (VCU) over multiple Controller Area Network (CAN) buses. Vehicle location is tracked by a time-resolved GNSS module, and an Intermodal Transport Control System (ITCS) provides the corresponding passenger counts and stop schedules. Data from these three sources (VCU, GNSS, and ITCS) are then synchronized and resampled. The second source is local meteorological data for the DTB service area from Open-Meteo, which draws on the ECMWF Integrated Forecasting System (ECMWF-IFS) model and assimilates observations from weather stations, aircraft, buoys, radar, and satellites to characterize the state of the atmosphere. The ECMWF-IFS model produces reliable short-term forecasts for the current day and the next few hours. A local edge data-fusion module performs timestamp alignment and data compression and transmits both streams to the cloud over 4G/5G. The two datasets are ingested offline, keyed on year-month-day-hour, aligned by timestamp, merged via a left join, and completed by interpolation (multi-source data aggregation). Continuous trajectories are then segmented into discrete "origin-stop-destination-stop" inter-stop trips. Finally, data cleaning, anomaly detection, missing-data imputation, and quality control are carried out sequentially, together with feature engineering, to produce the final modeling dataset.

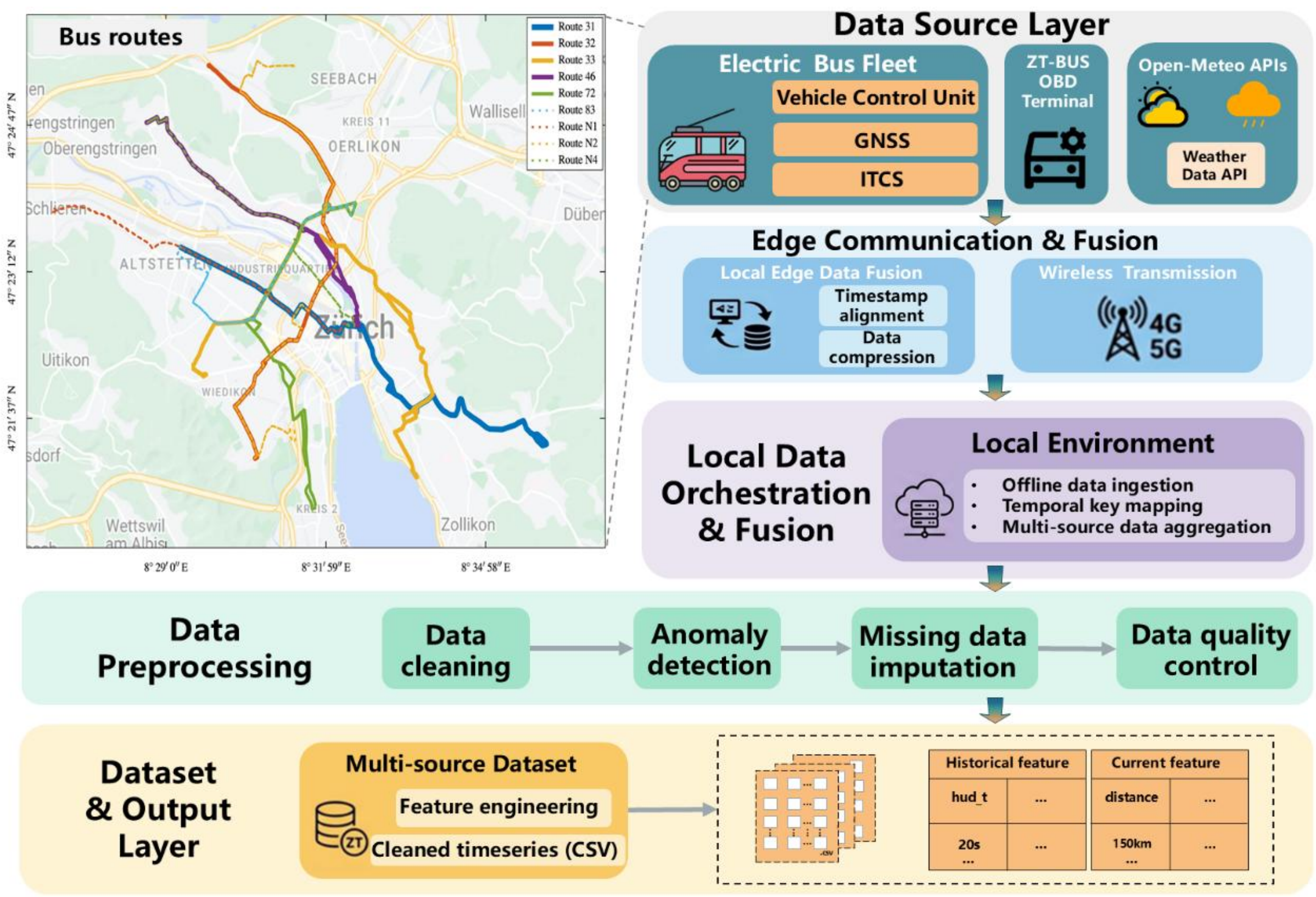

Fig. 1. DTB data collection and preprocessing workflow.

## 3.1 Data sources

The multi-dimensional dataset used in this study contains vehicle operation data and meteorological data.

### 3.1.1 Vehicle operation data (ZTBus dataset)

The base data come from the Zurich Transit Bus (ZTBus) dataset of the Zurich public-transport system. Two HESS lighTram® 19 DC dual-source trolleybuses (vehicle IDs 183 and 208) are equipped with dual power systems and can draw power either from the overhead catenary or from the on-board traction battery. The two DTBs run on lines 31, 32, 33, 46, 72, 83, N1, N2, and N4 of the Zurich transit network. The data cover numerous full-day driving missions between 2019 and 2022, with vehicle speed, traction force, distance, altitude, ambient temperature, door status, and total electric power demand recorded at 1 Hz by the VCU and Global Navigation Satellite System (GNSS). Passenger counts are estimated by an on-board infrared system. Each mission lasts about 8-10 hours and spans the whole operating day, providing sufficient temporal coverage to capture energy dynamics at morning, evening, and off-peak periods. The average driving speed is around 15 km/h, reflecting the frequent stops and signal interference typical of Zurich's inner city. The recorded average power demand ranges from about 15 to 35 kW, and the overall energy use lies between 1.5 and 2.0 kWh/km, which matches physical expectations for urban transit. Fig. 2 shows an example time-series segment from bus 183 on line VBZ 72 during a morning mission on 10 November 2022.

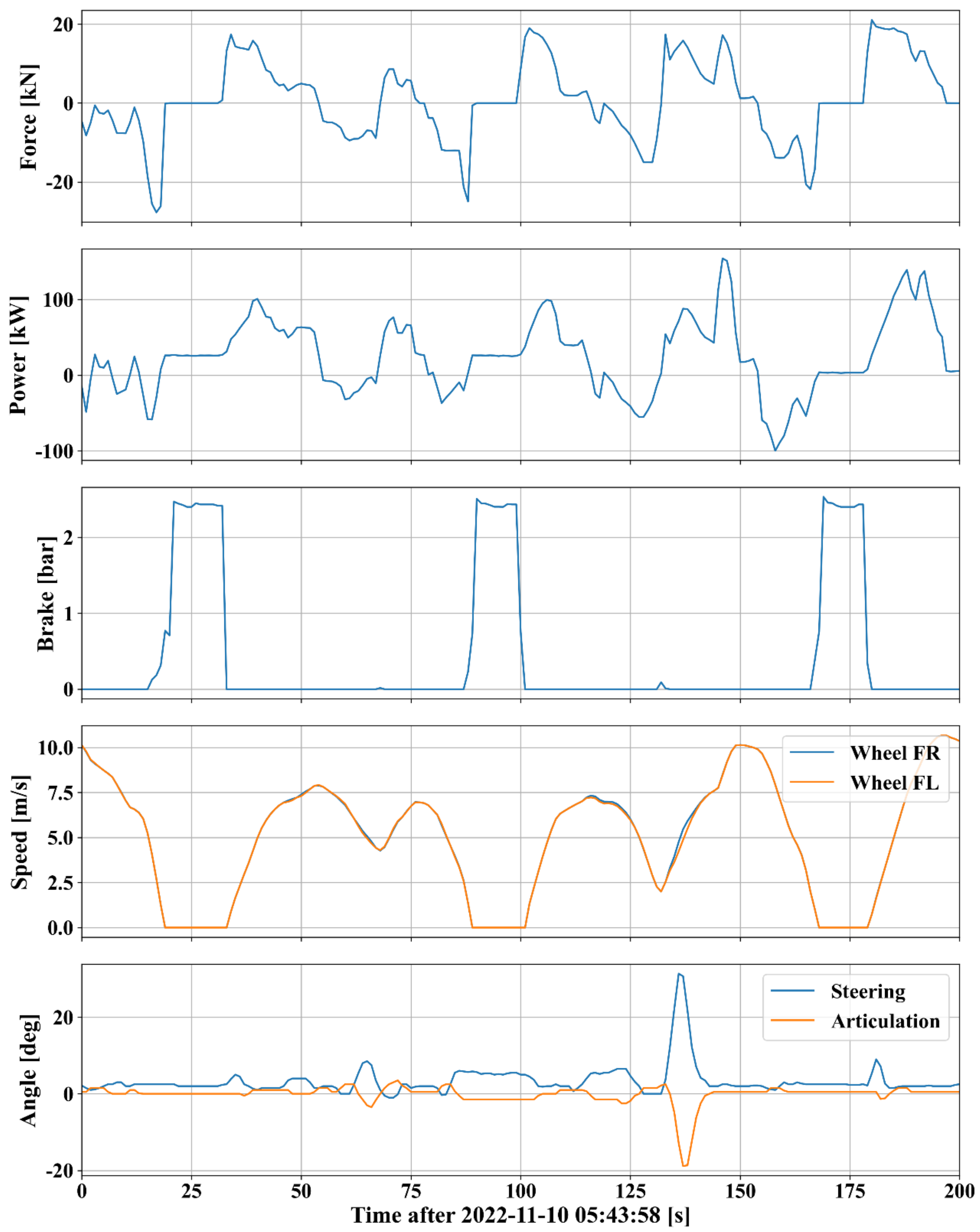

Fig. 2. Time-series view of bus 183 on line VBZ 72 on 10 November 2022.

### 3.1.2 External meteorological data

Because external conditions drive nonlinear fluctuations in both traction and auxiliary energy, this study adds high-resolution meteorological features that are matched to the vehicle trajectories in both space and time. Meteorological data are retrieved from the Open-Meteo Historical Weather API, which draws on the ERA5 reanalysis dataset of the European Centre for Medium-Range Weather Forecasts (ECMWF). Hourly records for the Zurich area are extracted to align with ZTBus operating periods, covering temperature, snow depth, precipitation, relative humidity, atmospheric pressure, and wind speed. These variables supplement environmental context that on-board sensors cannot provide.

## 3.2 Preprocessing and fusion

The vehicle operation data and the meteorological data are heterogeneous and are sampled at different rates (seconds vs. hours). A rigorous fusion and preprocessing pipeline is designed to ensure spatio-temporal consistency and statistical validity of the model inputs.

### 3.2.1 Spatio-temporal alignment and smooth fusion

To reconcile the sampling-frequency gap, year, month, day, and hour are extracted from the vehicle timestamps as the alignment key, and the vehicle data are rolled up to hourly granularity. Meteorological features are joined to the vehicle records via a left join, so that valuable vehicle samples are not discarded due to gaps in the meteorological series. Residual gaps are handled by linear interpolation together with forward and backward filling, preserving the integrity of the vehicle trajectories while smoothing the meteorological time series.

### 3.2.2 Target reconstruction and ratio-based features

To remove the effect of section length on the absolute energy value, the prediction target is set to the standardized energy per kilometer (Energy_per_km), defined as the total section energy divided by the distance traveled. Absolute values of time, distance, and mechanical work are further converted into relative ratios—for example, the share of each phase in total time and the fraction of traction work in total mechanical work. This normalization decouples the features from route-specific physical constraints, allowing the model to focus on driving-behavior patterns as determinants of energy use and to generalize to unseen routes.

### 3.2.3 Leakage-free standardization and temporal slicing

To accelerate convergence and remove unit effects, all continuous features are standardized. Because the target is right-skewed and long-tailed, a logarithm of one plus x (log1p) transformation is applied before Z-score standardization to reduce the impact of extreme outliers on gradient updates. Standardization statistics are computed from the training set only and reused on the validation and test sets to prevent leakage. For features that reflect intrinsic properties of specific segments, historical averages for the same origin-destination pair are used as calibration. On this basis, feature matrices that contain both "previous-step" and "current-step" information are constructed as structured inputs for the time-aware model.

The fused inter-stop dataset contains eight feature categories — time, operation, route, traction phase, coasting phase, braking phase, weather, and energy — totaling 63 fields (see Appendix Table A2).

# 4. Methodology

## 4.1 Problem statement

Inter-stop energy prediction for a DTB is the task of estimating the energy used over a specific between-stop section from multi-dimensional vehicle state data. The full process from one stop to the next is defined as an "inter-stop trip", and its energy prediction is essentially a supervised regression problem that fuses static attributes with dynamic temporal features [32]. Given the current time step t, the energy of the current trip is predicted from the states of the previous u trips together with the static tabular and dynamic temporal features of the current trip:

$$\hat{y}_t = f_\theta\left(x_t^{\mathrm{tab}}, x_t^{\mathrm{seq}}, y_{t-u,t-1}\right) \tag{1}$$

where $\hat{y}_t$ is the predicted energy of the t-th section; $f_\theta$ is the nonlinear mapping parameterized by θ (that is, the Time-Aware TabM model); $x_t^{\mathrm{tab}}$ denotes the static tabular features of the current trip; $x_t^{\mathrm{seq}}$ denotes the dynamic sequential features of the current trip; and $y_{t-u,t-1}$ is the historical memory from the past u trips.

## 4.2 Overall framework

The proposed framework for DTB energy prediction and causal analysis consists of three modules (see Fig. 3): (1) multi-source data fusion and preprocessing, which aligns ZTBus telemetry with Open-Meteo meteorological records in space and time, constructs the heterogeneous dataset, and performs feature engineering; (2) time-aware prediction and ensemble modeling, which uses TabM as the backbone, embeds Time2Vec periodic encoding, adopts a flattened GEMM batch-ensemble, and tunes hyperparameters by Bayesian optimization; and (3) interpretability and causal inference, which uses SHAP for attribution, DirectLiNGAM for causal-structure discovery, and a T-Learner for average-treatment-effect estimation. The core components are described below.

### 4.2.1 Construction of the heterogeneous dataset

**(1) Spatio-temporal alignment.**

Vehicle telemetry is sampled at the second level, while meteorological data are provided at the hourly level. A timestamp-indexed alignment strategy is used. Vehicle records are tagged with a structured time key $(y, m, d, h)$, and meteorological records are aggregated at the same hourly resolution:

$$\bar{x}_{\mathrm{clima}}(y, m, d, h) = \frac{1}{\left|S_{\mathrm{y,m,d,h}}\right|} \sum_{t \in S_{\mathrm{y,m,d,h}}} x_{\mathrm{clim}}(t) \tag{2}$$

where $\bar{x}_{\mathrm{clima}}(y, m, d, h)$ is the hourly mean meteorological value and $S_{\mathrm{y,m,d,h}}$ is the set of meteorological observations in that window. After a left join on the composite key $(y, m, d, h)$, residual gaps are filled by linear interpolation combined with forward and backward filling.

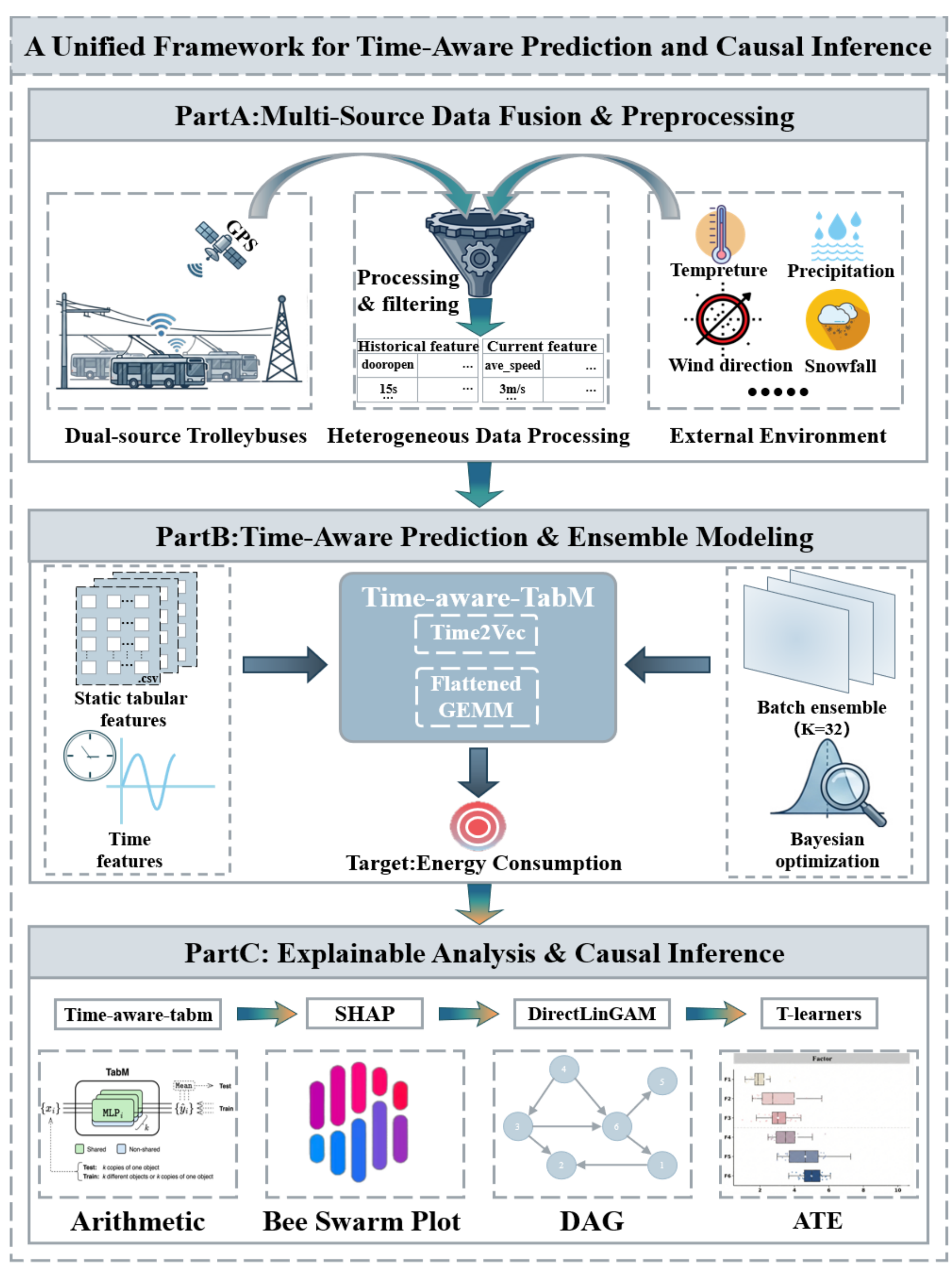

Fig. 3. Unified analysis framework for DTB energy drivers.

**(2) Feature engineering.**

Static attributes are one-hot encoded, continuous variables are Z-score standardized, and Time2Vec is used to map timestamps into vectors that carry both periodic and trend components. Because the raw variables span multiple physical domains — time, distance, speed, and energy — with markedly different scales and distributions, all continuous features are standardized as

$$\tilde{x} = \frac{x - u}{\sigma} \quad (3)$$

where $u$ and $\sigma$ are the mean and standard deviation estimated on the training partition only and reused for the validation and test sets to avoid leakage, and $\tilde{x}$ is the standardized value.

Because the target variable (energy per km, kJ/km) is right-skewed, a logarithmic transformation $x' = \ln(1 + x)$ is applied before standardization to compress the long tail:

$$\tilde{x}_{\text{energy}} = \frac{\ln\left(1 + x_{\text{energy}}\right) - u_{\log}}{\sigma_{\log}} \quad (4)$$

where $u_{\log}$ and $\sigma_{\log}$ are the mean and standard deviation on the training set and $\tilde{x}_{\text{energy}}$ is the standardized output.

### 4.2.2 Time-Aware TabM prediction model

**(1) TabM backbone.**

TabM is a parameter-efficient ensemble method for deep tabular modeling that combines batch ensembling with an MLP backbone [33]. By sharing backbone weights and introducing lightweight independent ensemble heads, TabM performs parallel inference of multiple sub-models in a single forward pass, capturing the diversity of the data distribution and reducing prediction variance at negligible extra cost. This weight-sharing mechanism mitigates the overfitting that typically affects single neural networks on tabular data while preserving inference efficiency.

**(2) Time2Vec embedding.**

In transit energy prediction, time features carry important periodic information about traffic flow. Traditional one-hot encoding ignores the continuity and periodicity of time, while naive numerical normalization cannot express periodic patterns. Inspired by the Fourier transform, Time2Vec maps a scalar timestamp to a low-dimensional continuous vector that contains both a linear component for long-term trend and sinusoidal components for periodic patterns. For an input time feature τ, the i-th embedding dimension is

$$\tau[i] = \begin{cases} \omega_i \tau + \varphi_i, \text{if} i = 0 \\ \sin(\omega_i \tau + \varphi_i), \text{if} 1 \leq i \leq k \end{cases} \quad (5)$$

where $\omega_i$ and $\varphi_i$ are learnable frequency and phase parameters and $k$ is the number of periodic terms. The linear term $(i = 0)$ captures non-periodic trends such as battery aging, while the periodic terms $(1 \leq i \leq k)$ use sinusoidal activations to capture daily or weekly cycles such as the tidal pattern of morning and evening peaks. Through this module, the model learns time-dependent patterns at multiple frequencies and concatenates them with static tabular features, improving its sensitivity to the dynamic traffic environment.

**(3) Flattened-GEMM accelerated batch ensemble.**

TabM uses batch ensembling to improve generalization and to quantify predictive uncertainty. Unlike classical ensembles that train multiple independent sub-models and suffer from parameter explosion, this paper uses a flattened general matrix multiplication (GEMM) strategy to realize ensemble learning with very high parameter efficiency in a single forward pass.

Batch ensembling elegantly resolves the tension between parameter sharing and model diversity. It shares the main weight matrix $\boldsymbol{W}$ and, for each sub-model $k$, introduces lightweight rank-one scaling vectors ($\mathbf{r}_k$ and $\mathbf{s}_k$) to approximate an independent weight matrix $\boldsymbol{W}_k$:

$$\boldsymbol{W}_k = \boldsymbol{W} \circ (\mathbf{s}_k \mathbf{r}_k{}^{\mathrm{T}}) \tag{6}$$

where $\circ$ denotes element-wise multiplication. This design reduces the ensemble parameter complexity from $O(K \times d_{\mathrm{in}} \times d_{\mathrm{out}})$ to $O(d_{\mathrm{in}} \times d_{\mathrm{out}} + K(d_{\mathrm{in}} + d_{\mathrm{out}}))$, sharply lowering memory consumption while providing weight-sharing regularization.

To exploit GPU-friendly dense matrix operations and avoid the fragmentation of per-model loops, the ensemble dimension is folded into the batch dimension. Given a batch input $\mathbf{X} \in \mathbb{R}^{B \times D}$, it is first expanded to $\mathbf{X}' \in \mathbb{R}^{B \times K \times D}$, multiplied element-wise by the scaling factor $\mathbf{r}$, and reshaped into $\widehat{\mathbf{X}} \in \mathbb{R}^{(B \cdot K) \times D}$, so that multi-model inference is formally reduced to a single-model large-batch computation. A standard matrix multiplication $\mathbf{Y} = \widehat{\boldsymbol{X}} \boldsymbol{W}^{\mathrm{T}}$ is then performed, and the output is reshaped back to $\mathbb{R}^{B \times K \times D_{out}}$, multiplied by the output scaling factor $\mathbf{s}$, and averaged over the ensemble dimension to obtain the final prediction. This design achieves a representation capacity comparable to training $K$ independent deep networks with minimal extra parameters.

### 4.2.3 Bayesian hyperparameter optimization

The performance of Time-Aware TabM depends strongly on hyperparameter settings. Because the model involves multiple hyperparameters — network depth, width, regularization strength, and learning rate — grid search becomes expensive and random search lacks adaptivity [7]. A Bayesian optimization strategy based on the TPE is therefore adopted and implemented with Optuna.

Unlike standard Sequential Model-Based Optimization (SMBO), which directly models $p(y \mid x)$ [7], TPE uses Bayes' rule in reverse. A threshold $y^*$ (a quantile of the currently observed objective) is set, and the historical trials are split into "good" and "poor" groups, whose conditional densities $l(x)$ and $g(x)$ are modeled separately:

$$p(x|y) = \begin{cases} l(x), \text{if} y < y^{\Psi} \\ g(x), \text{if} y \geq y^{\Psi} \end{cases} \tag{7}$$

where $l(x)$ is the parameter density for trials whose objective is below the threshold (better-performing) and $g(x)$ is the density for worse-performing trials.

TPE then chooses the next configuration by maximizing the expected improvement, which is equivalent to maximizing the ratio $l(x)/g(x)$, allowing the search to converge efficiently to a neighborhood of the global optimum. A median pruner is used at the same time: validation Mean Absolute Percentage Error (MAPE) is monitored during training, and trials whose performance falls below the historical median are terminated early, so that compute resources are concentrated on the more promising configurations. The overall architecture is shown in Fig. 4.

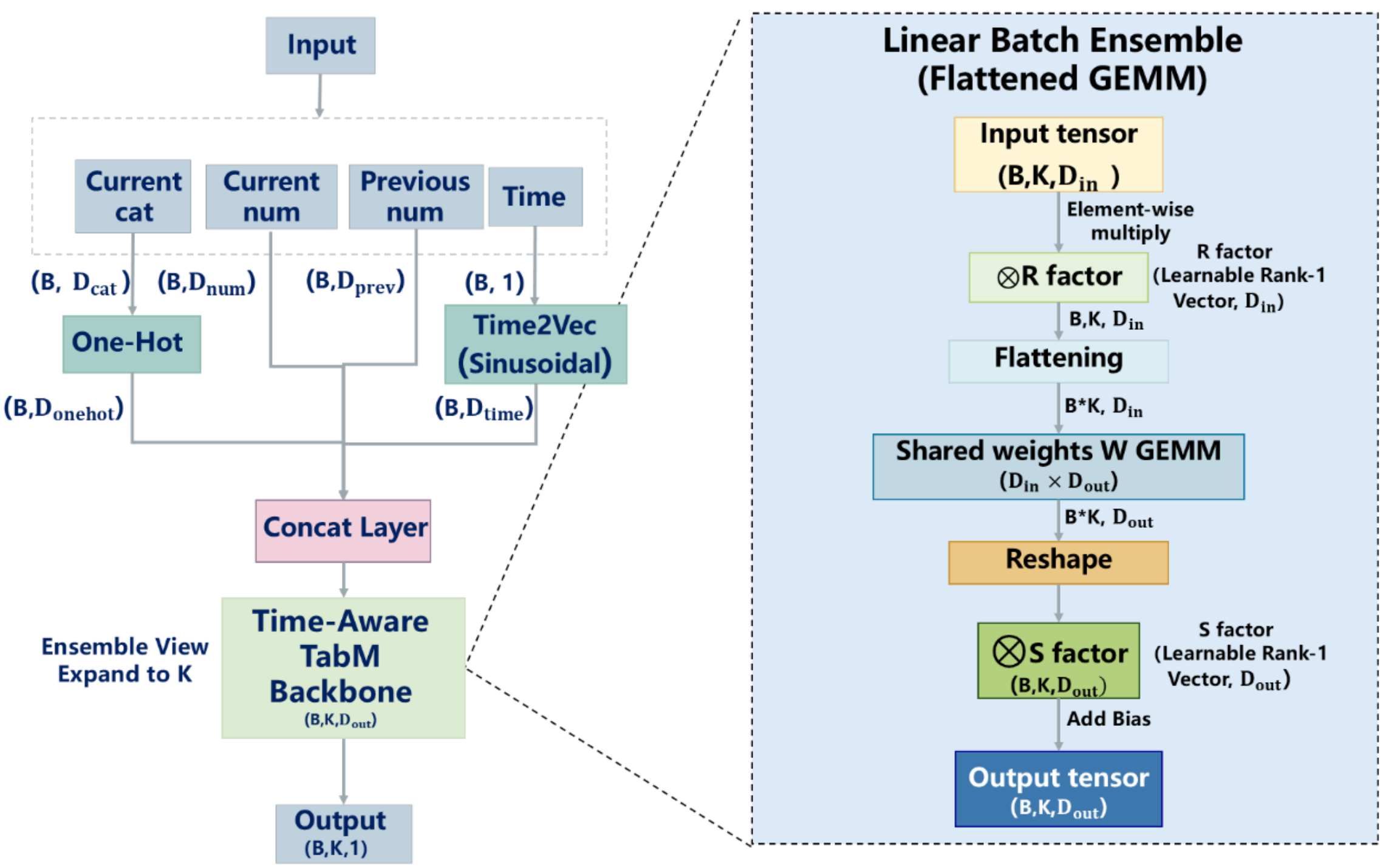

Fig. 4. Model architecture.

### 4.2.4 Causal inference

Building on the high-accuracy predictive model, a causal-inference framework is introduced to uncover the drivers of energy use from both the structural-discovery and the effect-estimation perspective.

**(1) Causal structure discovery via DirectLiNGAM.** DirectLiNGAM [34] is used to build a directed acyclic graph (DAG) over the features and the energy variable. The algorithm is based on a linear non-Gaussian acyclic model and exploits the non-Gaussianity of the observed data to identify causal directions without interventional experiments. Let $\mathbf{x} = [x_1, x_2, \cdots, x_m]^{\mathrm{T}}$ be the standardized variables; the generating process is

$$\mathbf{x} = \boldsymbol{B}x + \mathbf{e} \tag{8}$$

where $\boldsymbol{B}$ is an $m \times m$ adjacency matrix whose element $b_{ij}$ represents the direct causal effect of $x_j$ on $x_i$, and $\mathbf{e} = [e_1, \cdots, e_m]^T$ is a vector of mutually independent, non-Gaussian exogenous noises. Under the correct causal ordering, $\boldsymbol{B}$ can be permuted into a strictly lower-triangular matrix ($\boldsymbol{b}_{ii} = 0$), ensuring the acyclicity assumption.

DirectLiNGAM determines the causal order through successive regression and independence testing. For any two variables $x_i$ and $x_j$ whose direction is not yet determined, a linear regression is performed,

$$x_j = c_j^{(i)} x_i + r_j^{(i)} \tag{9}$$

where $c_j^{(i)} = \mathrm{cov}(x_i, x_j)/\mathrm{var}(x_i)$ and $r_j^{(i)}$ is the residual. The key criterion is that if $x_i$ is a cause of $x_j$, then the residual $r_j^{(i)}$ is statistically independent of $x_i$.

At each iteration, the algorithm uses mutual information or a kernel-based independence criterion to find the root variable whose regression residuals on the remaining variables are the most independent, places it at the top of the causal chain, and subtracts its effect from the remaining variables:

$$x^{(s+1)} = x^{(s)} - c^{(root)} x_{root}^{(s)} \tag{10}$$

where $s$ and $s$+1 index the iteration. This recursion continues until all variables are ordered, so that the causal topology of energy use is reconstructed adaptively.

**(2) Treatment-effect estimation via T-Learner.** Once the causal structure is established, a T-Learner meta-learning framework [10] is used to quantify the net causal effect of each driver on energy use. Continuous features $V$ are typically non-uniformly distributed, so the 75th percentile is used as a threshold to discretize them into a binary treatment variable $T$:

$$T_i = \begin{cases} 1, \text{if} V_i < Q_{0.75}(V) \\ 0, \text{otherwise} \end{cases} \tag{11}$$

The T-Learner completely separates the treatment group (T = 1) from the control group (T = 0) and fits two independent base learners. Here, Time-Aware TabM is used as the base model, and the two learners estimate the conditional expectation given covariates $x$ (all other features except the treatment):

$$u_0(x) = E[Y|\mathrm{T} = 0|X = x] \tag{12}$$

$$u_1(x) = E[Y|\mathrm{T} = 1|X = x] \tag{13}$$

where $X = x$ represents covariates, which are state characteristics other than the intervention variable. The model $u_0$ is trained only on the control data (T = 0), and $u_1$ only on the treatment data (T = 1). At inference time, the difference between the two predictions gives the conditional average treatment effect (*CATE*) for sample $i$ with covariates $\mathbf{x}_i$:

$$\widehat{CATE}(x) = \hat{u}_1(x) - \hat{u}_0(x) \tag{14}$$

Averaging over the $N$ test samples yields the average treatment effect (*ATE*):

$$\widehat{ATE} = \frac{1}{N}\sum_{i=1}^{n} \widehat{CATE}(x) \tag{15}$$

where $N$ is the total sample size. $\widehat{ATE}$ reflects the average change in energy use caused by the treatment under the current covariate distribution. A positive $\widehat{ATE}$ indicates that the factor increases overall energy use, while a negative value indicates an energy-saving effect.

# 5. Experimental Analysis

## 5.1 Experimental setup

**(1) Software and hardware.** Experiments were run on a Windows 10 mobile workstation with an AMD 10-core CPU, 16 GB RAM, and an NVIDIA GeForce RTX 4060 Laptop GPU (8 GB), using Python 3.x and PyTorch with CUDA 12.8.

**(2) Dataset partition.** To evaluate generalization under different operating conditions, the data are split into training, validation, and test sets in a 60:20:20 ratio by random sampling. Random splitting preserves the joint distribution over season, temperature, and traffic flow across the three subsets and avoids evaluation bias from period concentration. The training set is used for gradient updates, the validation set for hyperparameter tuning and early stopping, and the test set only for final evaluation.

**(3) Model settings.** The model is trained with AdamW and a weight decay of 1e-5. Because the hyperparameter search space mixes continuous and discrete variables, grid search is infeasible; a TPE-based Bayesian optimization with median pruning is therefore used to terminate poorly converging trials. The search space is summarized in Table 1.

Table 1. Hyperparameter search space.

| No. | Parameter | Symbol | Range | Role |
|---|---|---|---|---|
| 1 | Ensemble size | k | 32 (fixed) | Number of TabM ensemble members |
| 2 | Learning rate | lr | 1e-4 - 1e-2 (log) | Weight update step size |
| 3 | Batch size | B | 128, 256, 512 | Samples per training batch |
| 4 | Number of blocks | n_blocks | 1, 2, 3, 4 | Model depth |
| 5 | Hidden dimension | d_block | 64, 128, 256, 512 | Hidden-layer width |
| 6 | Dropout rate | dropout | 0.0 - 0.5 | Regularization against overfitting |
| 7 | Weight decay | wd | 1e-6 - 1e-3 (log) | L2 regularization strength |

The Bayesian optimization is run on the validation set with a mixed search space over learning rate, batch size, depth, and hidden width. The experiment uses 50 trials, each trained for at most 50 epochs, with early stopping (patience = 15) and median pruning. The best configuration after 50 trials is: learning rate 0.00536, batch size 128, 4 network blocks, hidden dimension 256, and dropout 0.04.

## 5.2 Baseline models

To evaluate Time-Aware TabM comprehensively, ten representative algorithms are selected as baselines, covering statistical models, classical machine learning, and deep learning.

1) Statistical model: ① ARIMA [35] is a classical linear time-series baseline that fits trend and seasonality through autoregression, differencing, and moving averages.

2) Classical machine learning: ② SVR [36] uses kernel functions for high-dimensional nonlinear mapping;

③ RF [37], ④ XGBoost [38], and ⑤ LightGBM [39] represent bagging and boosting ensembles, with LightGBM using Gradient-based One-Side Sampling (GOSS) and Exclusive Feature Bundling (EFB) to improve efficiency on large data; ⑥ DLN [40] is a stacked-autoencoder plus MLP feed-forward network.

3) Deep learning: ⑦ a 1-D residual CNN [41] captures local temporal features through convolution; ⑧ LSTM [42] uses gating to model long-range dependencies; ⑨ TFT [32] is a Transformer-based architecture for temporal feature interaction; ⑩ TabPFN [43] is a Transformer-based tabular foundation model that performs fast inference on small datasets through Bayesian meta-learning. Together, these baselines cover linear and nonlinear, shallow and deep modeling paradigms, providing a systematic reference for evaluation.

## 5.3 Evaluation metrics

Mean Absolute Error (MAE), Root Mean Squared Error (RMSE), MAPE, the coefficient of determination $R^2$, and computation time (s) are used as evaluation metrics, with MAPE as the core optimization metric measuring relative accuracy across operating conditions. Definitions are given in Eqs. (16) - (19).

$$MAE = \frac{1}{n}\sum_{i=1}^{n}|y_i - \hat{y}_i| \tag{16}$$

$$RMSE = \sqrt{\frac{1}{n}\sum_{i=1}^{n}(y_i - \hat{y}_i)^2} \tag{17}$$

$$MAPE = \frac{100\%}{n}\sum_{i=1}^{n}\left|\frac{y_i - \hat{y}_i}{y_i}\right| \tag{18}$$

$$R^2 = 1 - \frac{\Sigma_{i=1}^{n}(y_i - \hat{y}_i)^2}{\Sigma_{i=1}^{n}(y_i - \bar{y})^2} \tag{19}$$

where $\hat{y}_i$ is the predicted inter-stop energy (kJ) for sample $i$,$y_i$ is the ground truth (kJ), and $n$ is the sample size.

## 5.4 Ablation study

Four variants are designed to evaluate the contribution of each core component: TabM-v0 (plain TabM with no enhancements), TabM-v1 (flattened GEMM only), TabM-v2 (Time2Vec and historical features only), and TabM-v3 (the full Time-Aware TabM). All variants use the hyperparameters obtained by Bayesian optimization and the same data split. Results are reported in Table 2 (CT denotes computation time).

Overall, the components contribute to the model in clearly different ways: the time-awareness module is the main driver of accuracy, while the compute-layer optimization only delivers its potential once the feature representation is strong enough. Looking at the compute-layer optimization, adding flattened GEMM alone (v1) actually degrades the model slightly relative to the baseline (v0): MAE rises from 403.69 to 407.77, MAPE from 7.68% to 7.75%, and $R^2$ drops from 0.975 to 0.974. This confirms that, without effective feature representation, pure matrix-operation acceleration does not translate into accuracy gains — faster arithmetic and stronger modeling capacity are not causally linked.

Table 2. Ablation results.

| Model | MAE | RMSE | MAPE | MAPE_std | MAPE_range | R² | CT |
|---|---|---|---|---|---|---|---|
| TabM-v0 | 403.69 | 536.60 | 7.68% | 0.12 | 0.58 | 0.975 | 10.18 s |
| TabM-v1 | 407.77 | 542.03 | 7.75% | 0.14 | 0.72 | 0.974 | 9.78 s |
| TabM-v2 | 353.43 | 472.16 | 6.68% | 0.08 | 0.40 | 0.980 | 18.80 s |
| TabM-v3 | 340.42 | 457.02 | 6.52% | 0.07 | 0.49 | 0.982 | 28.91 s |

The time-awareness module brings the largest leap. Fusing Time2Vec with historical state features (v2) cuts MAE to 353.43 (−12.45%) and MAPE to 6.68% (−1.00 %), and pushes $R^2$ from 0.975 to 0.980. At the same time, the standard deviation of MAPE drops from 0.12 to 0.08, and its range from 0.58 to 0.40, showing a real gain in both accuracy and stability. These results confirm that effective capture of temporal periodicity and historical operating patterns is the principal driver of accuracy improvement. Time2Vec encodes time information into learnable periodic and non-periodic features, enabling the model to capture the intrinsic temporal rhythms of traffic flow—including the tidal dynamics of morning and evening peak periods.

From the module interaction perspective, v3 further compresses MAPE to 6.52% and lifts $R^2$ to 0.982, giving the best overall performance. Interestingly, although flattened GEMM alone is slightly worse than the baseline, when combined with Time2Vec it pushes v3 below the Time2Vec-only v2 by 0.16 % of MAPE, indicating a complementary gain. A plausible explanation is that the parameter-sharing regularization of flattened GEMM suppresses overfitting across the ensemble members once the feature representation is already rich enough, letting the diversity of ensemble predictions play out more fully. The MAPE standard deviation of v3 (0.07) is also the lowest among all variants, confirming that it has the most robust generalization. The inference time (28.91 s) is about 2.8 times the baseline, but this cost is acceptable for off-line energy analysis and scheduling given the clear gains in accuracy and stability.

## 5.5 Prediction performance comparison

To eliminate the impact of randomness on the evaluation, five base random seeds are preset and each setting is run 30 times independently (150 runs in total), and both accuracy and variability are reported. All baseline hyperparameters are tuned by Bayesian optimization. For the ensemble learners: RF uses n_estimators = 250, max_depth = 20, min_samples_split = 6; XGBoost uses n_estimators = 2500, learning_rate ≈ 0.0465, max_depth = 4; LightGBM uses n_estimators = 2600, learning_rate ≈ 0.0377, max_depth = 8, num_leaves = 42. SVR uses an RBF kernel with C ≈ 5.50 and epsilon ≈ 0.017. For the deep baselines: DLN has three hidden layers (512-256-128); CNN uses hidden_dim = 64 and batch_size = 256; LSTM uses num_layers = 1 and hidden_dim = 64; TFT uses d_model = 128, nhead = 8, and num_layers = 2. Deterministic models (SVR, TabPFN, ARIMA, etc.) are not affected by the random seed and their standard deviation and range are marked "-". Results are summarized in Table 3 (CT denotes computation time).

Overall, Time-Aware TabM is the best on every metric, with a MAPE of 6.52% and $R^2$ of 0.982, clearly

outperforming all baselines. The results can be analyzed in three tiers.

Table 3. Performance comparison across models.

| Model | MAE | RMSE | MAPE | MAPE_std | MAPE_range | $R^2$ | CT |
|---|---|---|---|---|---|---|---|
| ARIMA | 915.12 | 1332.27 | 16.16% | - | - | 0.843 | 12.59 s |
| SVR | 607.80 | 878.55 | 10.96% | - | - | 0.932 | 134.25 s |
| RF | 543.80 | 763.07 | 10.37% | 0.02 | 0.12 | 0.948 | 33.46 s |
| DLN | 814.13 | 1112.83 | 14.68% | 1.01 | 5.48 | 0.890 | 0.47 s |
| TabPFN | 992.71 | 1332.46 | 18.45% | - | - | 0.843 | 9.73 s |
| LightGBM | 415.65 | 567.64 | 7.89% | - | - | 0.971 | 4.99 s |
| XGBoost | 420.46 | 572.00 | 7.89% | - | - | 0.971 | 18.42 s |
| CNN | 426.12 | 572.26 | 7.90% | 0.29 | 1.68 | 0.971 | 14.87 s |
| LSTM | 404.53 | 541.90 | 7.59% | 0.12 | 0.64 | 0.974 | 4.65 s |
| TFT | 466.99 | 636.06 | 8.30% | 0.24 | 1.07 | 0.964 | 10.04 s |
| **TabM-v3** | **340.42** | **457.02** | **6.52%** | **0.07** | **0.49** | **0.982** | **28.91s** |

The first tier contains the weakest models (MAPE > 10%). ARIMA, at 16.16%, lacks the capacity to capture the nonlinear coupling between driving behavior and environmental factors — its linear autoregressive structure is simply not expressive enough. TabPFN (18.45%), as a pre-trained tabular foundation model, is not fine-tuned on this large, domain-specific dataset and therefore generalizes poorly. DLN (14.68%) has a too-simple architecture and unstable training (MAPE_std = 1.01), so its capacity cannot support effective modeling of high-dimensional features. SVR (10.96%) has nonlinear mapping ability, but kernel methods struggle to disentangle multi-factor interactions in high-dimensional space and are sensitive to hyperparameters. RF (10.37%) is an ensemble method, but its bagging strategy is weaker than gradient-boosted trees at modeling the multi-dimensional feature interactions in the DTB data.

The second tier consists of tree ensembles and deep models with MAPE around 7-9%. XGBoost and LightGBM both reach 7.89%, showing that gradient-boosted trees have an intrinsic advantage in capturing nonlinear relations across many features; their deterministic inference also gives a MAPE standard deviation of zero. CNN (7.90%) is comparable in accuracy, but its MAPE standard deviation (0.29) and range (1.68) are larger, suggesting that 1-D convolution is sensitive to initialization. LSTM (7.59%) outperforms CNN in temporal dependency modeling thanks to its gating, but its longitudinal information propagation ignores horizontal feature correlations between input variables. TFT (8.30%) has an attention mechanism but cannot fully exploit its sequence-modeling strengths on this tabular-dominant task.

The third tier is the proposed model. Time-Aware TabM reduces MAPE by 1.07 % relative to LSTM and by 1.37 % relative to XGBoost, and also has the lowest MAPE standard deviation (0.07) and range (0.49). This advantage comes from three complementary elements: Time2Vec encoding gives the model a sense of periodic traffic-flow structure, the TabM backbone handles heterogeneous tabular features efficiently. The flattened GEMM batch ensemble suppresses overfitting through parameter sharing. Together, these components capture nonlinear dependencies in both time and feature space, producing a clear lead in both accuracy and stability. Figs. 5, 6a, and 6b show the raincloud plot and training/test performance comparisons, respectively.

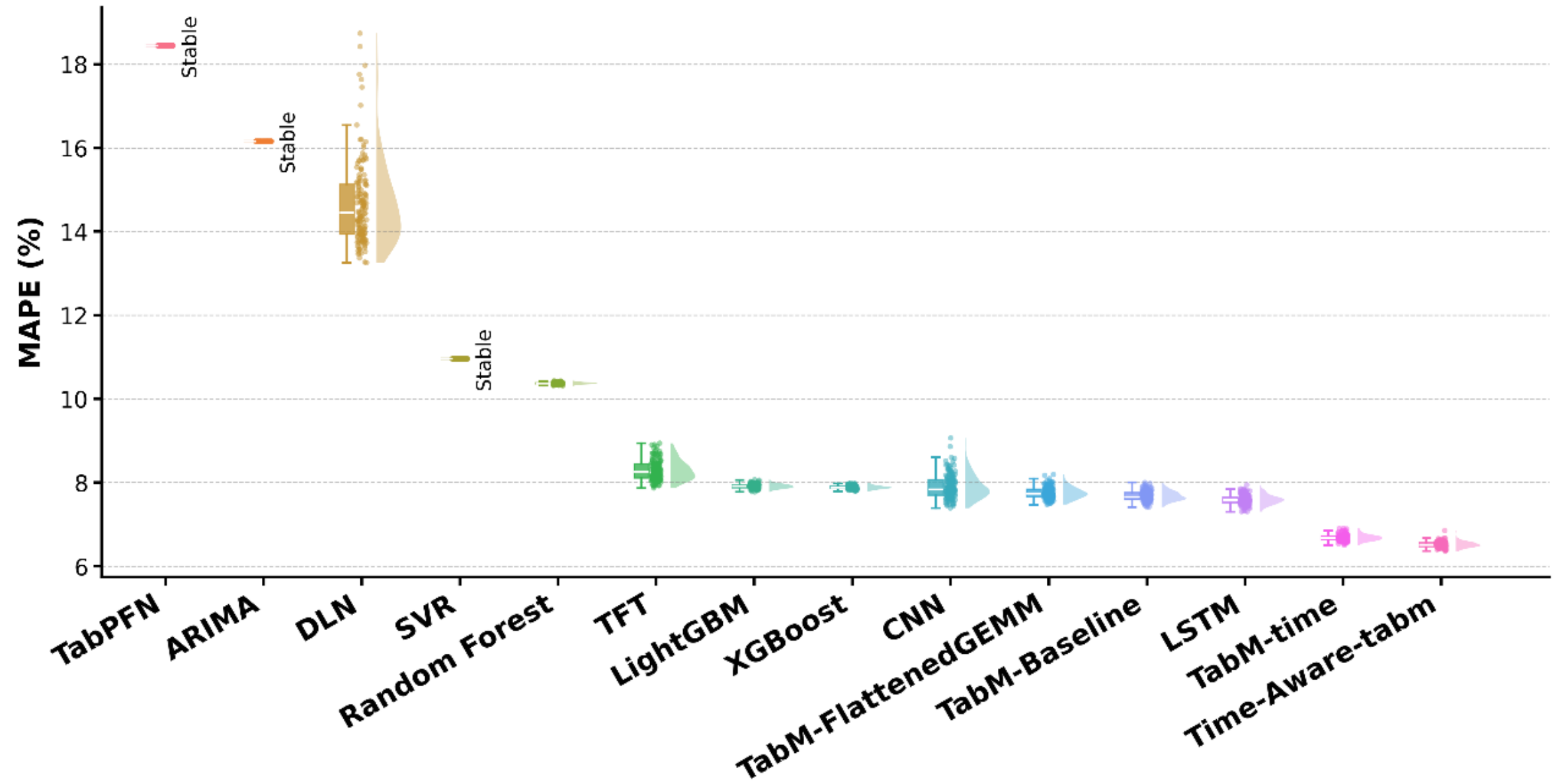

Fig. 5. Raincloud plot comparing models.

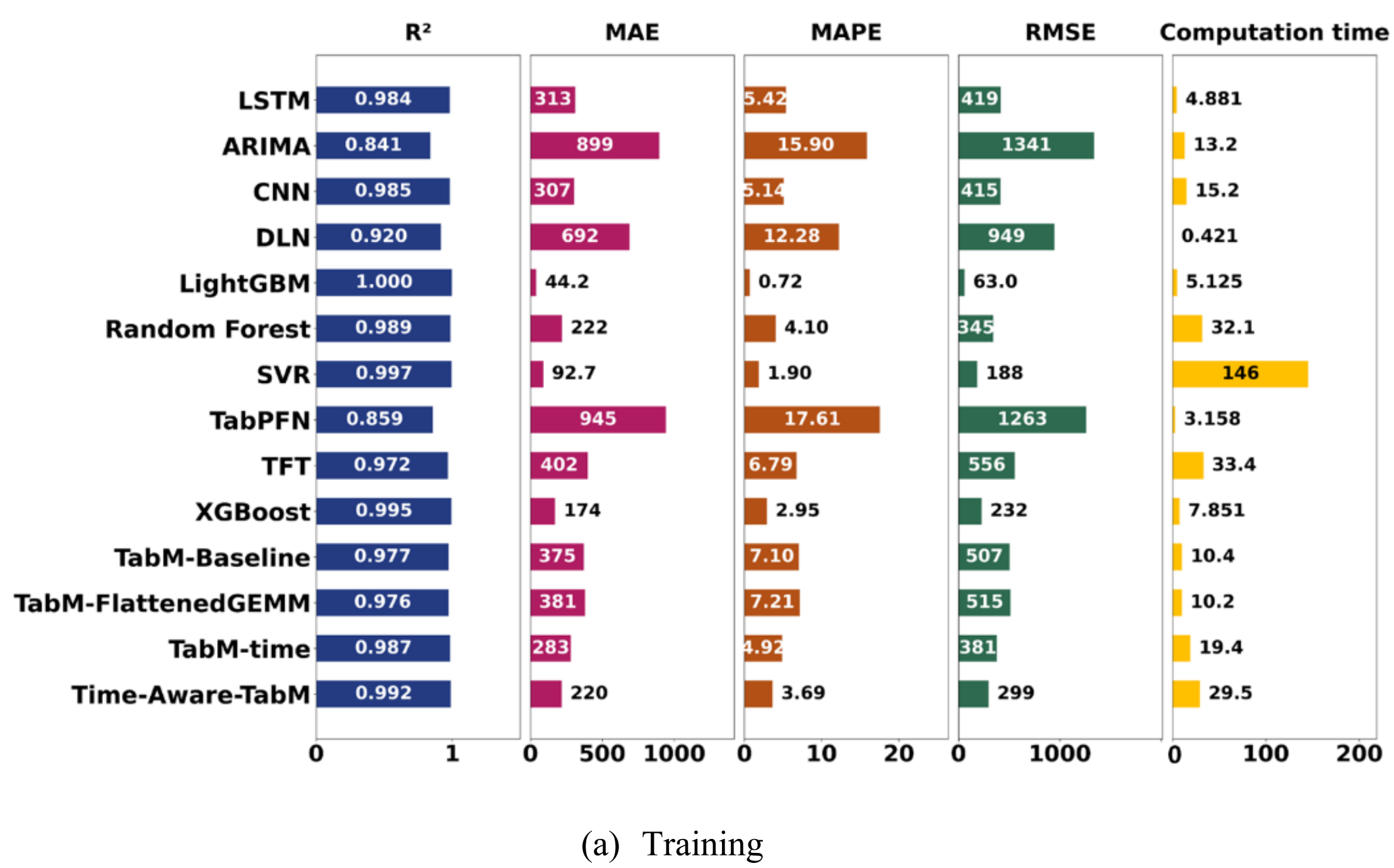

(a) Training

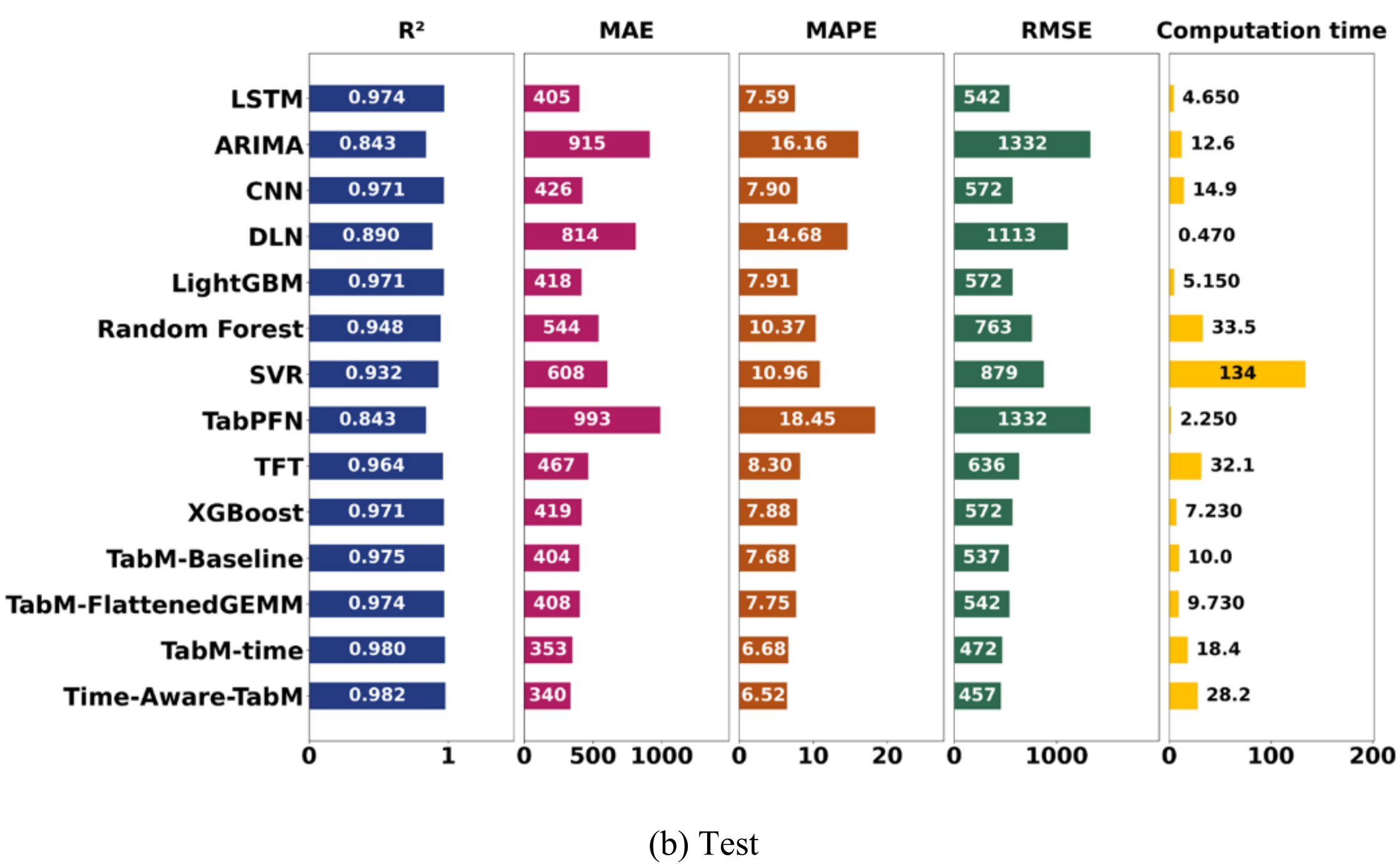

(b) Test

Fig. 6. Training and test performance comparison of all models. (a) Training; (b) Test.

## 5.6 SHAP interpretability analysis

To uncover the internal decision logic of Time-Aware TabM, SHAP-based attribution is performed on the test set and analyzed progressively at five levels: global importance, direction of impact, decision paths, individual attribution, and feature interactions.

### 5.6.1 Global feature importance (SHAP bar plot)

Fig. 7 shows the feature importance ranking based on the mean absolute SHAP value, together with a nested donut chart. In the inner ring of the donut, dark blue, light blue, and pink represent vehicle operating state, driving environment and conditions, and driving behavior, respectively; the outer ring corresponds to specific features whose colors match the bar chart. In terms of category contribution, vehicle operating state accounts for the largest share of total SHAP contribution (57.05%), followed by driving environment and conditions (25.87%) and driving behavior (17.08%). This distribution shows that the vehicle's own dynamic state is the most informative source for energy prediction.

At the feature level, the regenerative braking ratio (Faratio), average speed (Ave_speed), coasting time (Hua_t), dwell time (Stop_t), and coasting distance (Hua_d) rank in the top five. The contribution of the regenerative braking ratio is markedly larger than any other variable, reflecting the central role of kinetic energy recovery in DTB energy use. Average speed follows closely, in line with its role as a composite descriptor of driving conditions. Notably, both coasting time and dwell time enter the top five; both are tied to non-traction states, which highlights how energy losses during non-effective driving phases significantly shape overall energy efficiency.

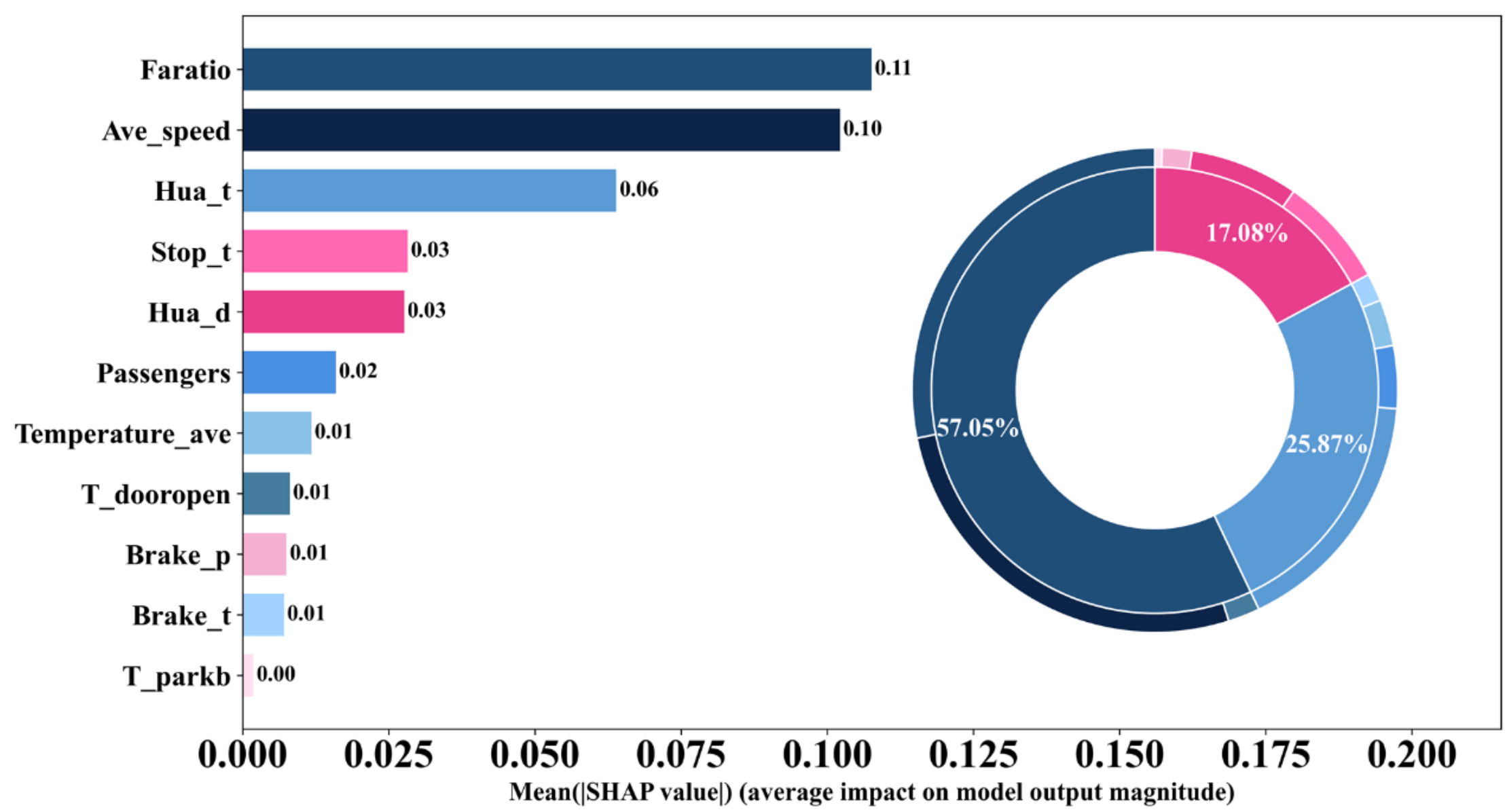

Fig. 7. Global feature importance.

### 5.6.2 Direction and distribution of feature impact (SHAP beeswarm plot)

The SHAP beeswarm plot in Fig. 8 reveals the mapping between feature values and the direction of energy change, showing four typical patterns. (1) The regenerative braking ratio (Faratio) and coasting time (Hua_t) exhibit the most pronounced negative SHAP contributions. Samples with high feature values concentrate in the negative SHAP region, indicating that a greater share of kinetic energy recovery and prolonged inertial coasting substantially reduce predicted energy consumption. The effect of Faratio is particularly consistent: its high-valued samples yield the largest negative SHAP magnitudes among all features. This aligns with the physical behavior observed in the ZTBus dataset, where braking power can reach approximately −320 kW, reflecting the regenerative braking system's capacity to convert kinetic energy into stored electrical energy. High-value coasting time samples also shift towards negative SHAP values, matching the mechanism in which inertial coasting reduces active traction effort. (2) Door-open time (T_dooropen) and coasting distance (Hua_d) show positive (energy-increasing) patterns: high values shift towards positive SHAP values. In non-traction states, longer door-open time keeps auxiliary systems — especially HVAC — running, which directly raises per-kilometer energy use. The positive effect of coasting distance indicates that overly long traction-free coasting is not an optimal energy-saving strategy; driver training should focus on precise timing of coasting to reduce wasted kinetic-energy dissipation. (3) Average speed shows a steep positive penalty in the low-speed range, which rapidly declines and flattens as speed increases, indicating that speed optimization yields the largest marginal energy savings under congestion. (4) Low ambient temperature (Temperature_ave) consistently shifts samples towards positive SHAP values, attributable to the additional auxiliary load from electric heaters under cold conditions. Passenger number (Passengers), somewhat counter-intuitively, shows a negative trend; the physical explanation is discussed in the causal-analysis section.

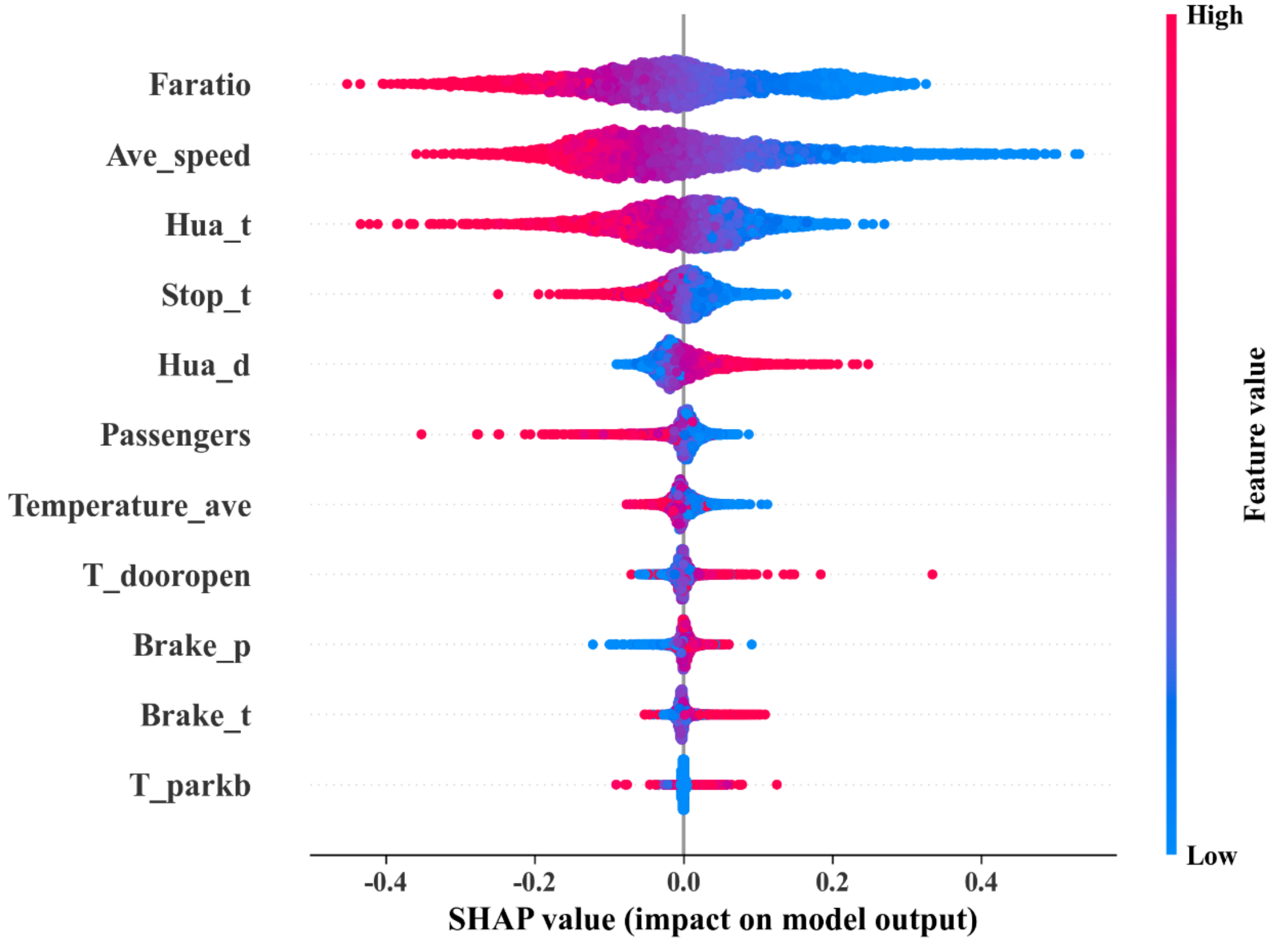

Fig. 8. SHAP beeswarm plot.

### 5.6.3 Prediction decision paths (SHAP decision plot)

The SHAP decision plot in Fig. 9 shows the full path along which the model accumulates feature contributions from the average baseline to the final prediction for the first twenty test samples. The paths show a clear hierarchical structure: at the bottom, for low-level operational features such as parking-brake time (T_parkb), brake time (Brake_t), and door-open time, all samples converge tightly, indicating only minor local perturbations. As the paths move up to mid-level features (ambient temperature and dwell time), the sample lines begin to diverge noticeably. At the top, for the dominant kinematic features (coasting distance and average speed), the paths fan out sharply: high-energy samples (purple-red) are pushed strongly in the positive direction, while low-energy samples (blue) receive a strong negative push.

This hierarchical structure delivers a key insight: the decisive split in energy prediction does not occur at the micro-operation level but is dominated by macro-level features related to kinematic state. This finding matches the global importance ranking in Section 5.6.1 and, from the path perspective, further confirms the dominant role of kinematic state in DTB energy prediction.

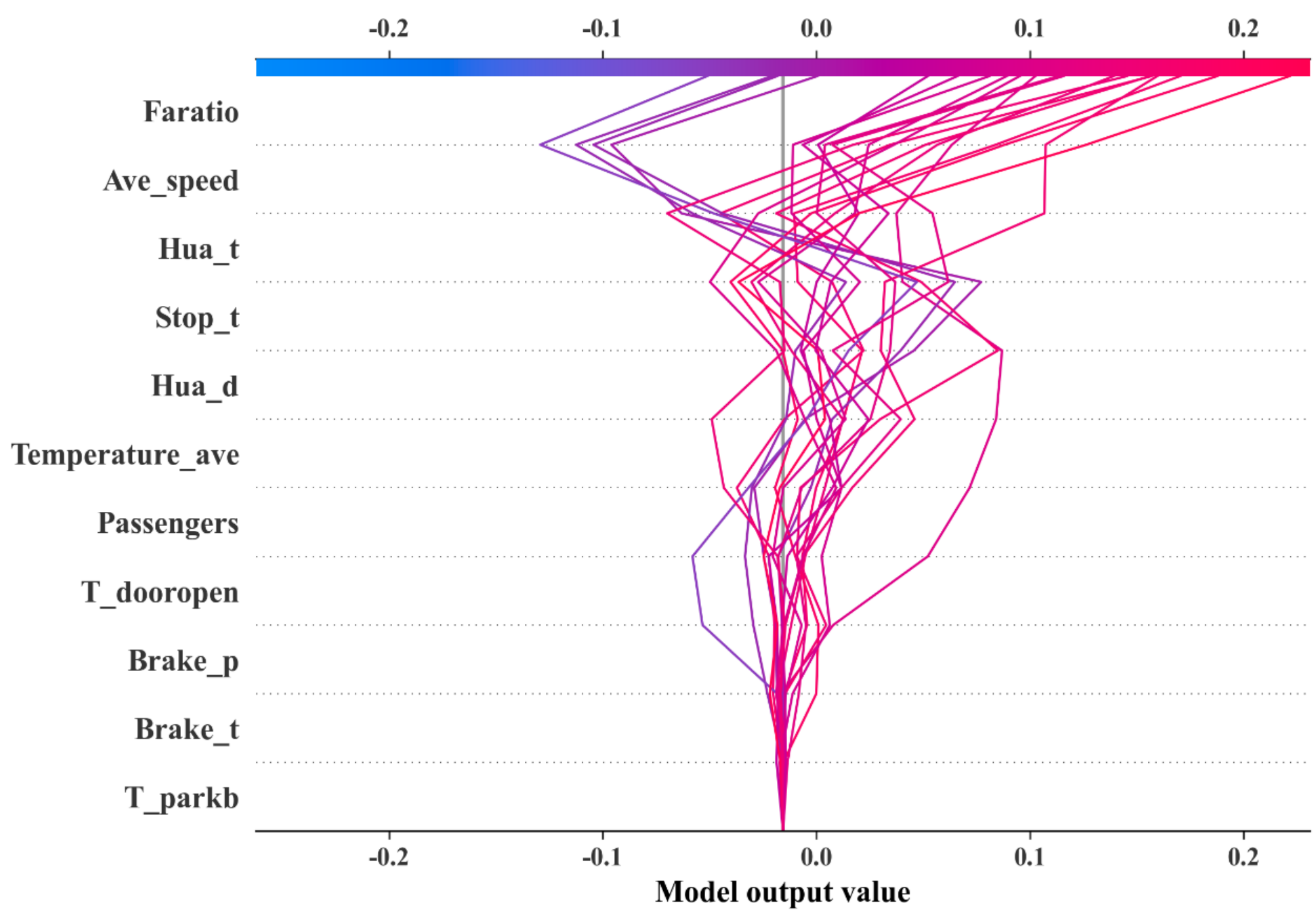

Fig. 9. SHAP decision plot.

### 5.6.4 Consistency of individual attribution

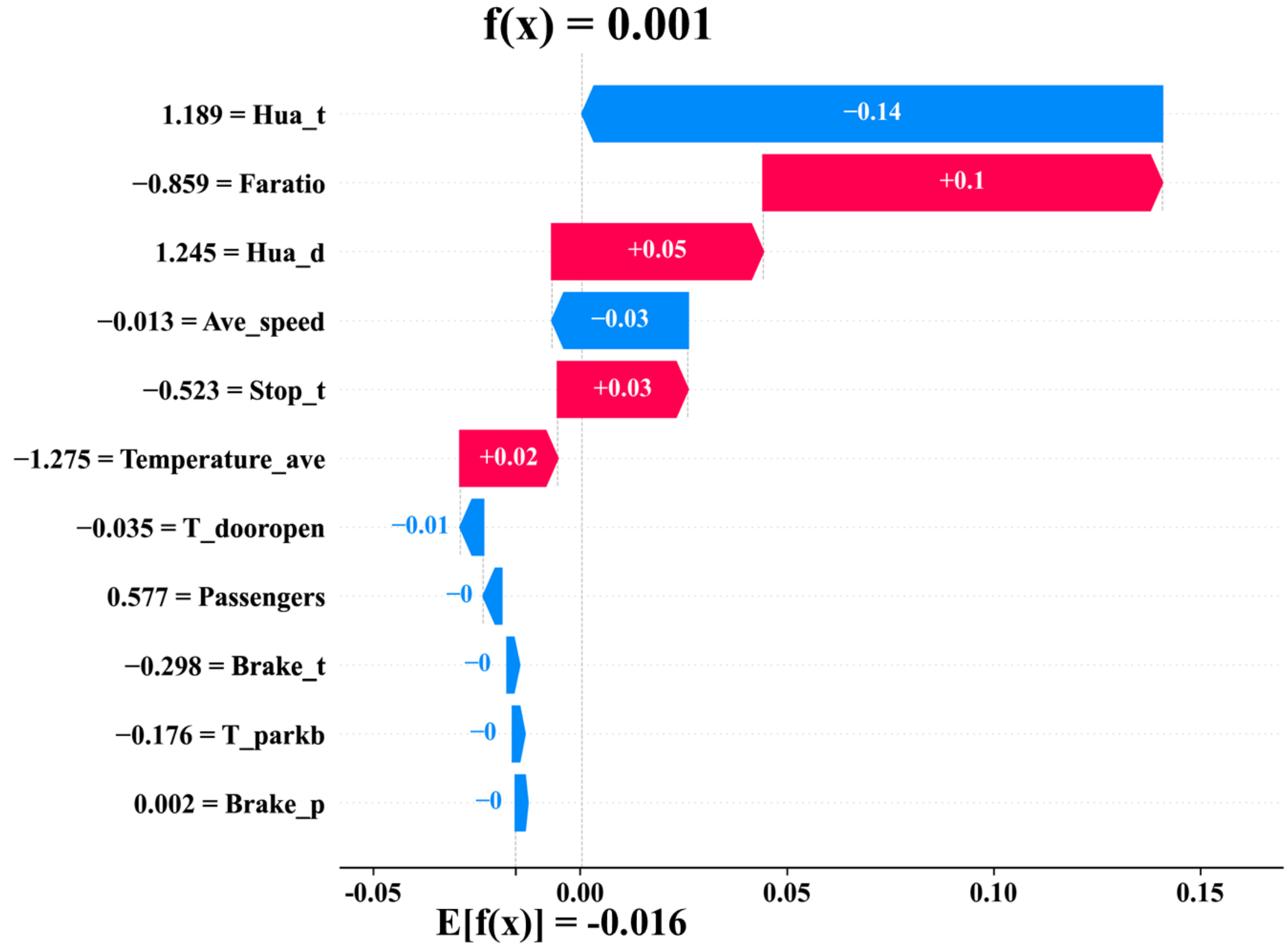

Fig. 10. SHAP waterfall plot.

To check whether the global patterns hold at the individual level, one sample is examined in detail. The SHAP waterfall plot (Fig. 10) shows that the final prediction (f(x) = 0.001) is slightly above the model baseline (E[f(X)] = −0.016). A detailed look at feature contributions shows that coasting time (Hua_t, value 1.189) contributes −0.14 in SHAP, making it the main driver pushing the prediction downward (energy saving). However, Faratio (value −0.859) and coasting distance (Hua_d, value 1.245) contribute +0.10 and +0.05 respectively, becoming the main forces pushing the prediction upward. Average speed (Ave_speed, value −0.013) then adds a small −0.03, while a shorter-than-average dwell time (Stop_t, value −0.523) and a lower average temperature (Temperature_ave, value −1.275) contribute +0.03 and +0.02, respectively. The interplay of these factors not only offsets the energy saving from the long coasting time but also nudges the final prediction slightly above the baseline.

The SHAP force plot (Fig. 11) confirms this analysis from a force-balance perspective. The numbers next to each feature are standardized input values; negative values indicate that the feature is below the training mean, and positive values indicate the opposite. In this sample, the "energy-saving camp" (blue) is led by the long coasting time (1.189), which directly cuts energy use; the "energy-increasing camp" (red) is led by the lower regenerative braking ratio (−0.8592) and longer coasting distance (1.2452), reinforced by a shorter dwell time (−0.5232) and lower average temperature (−1.2754).

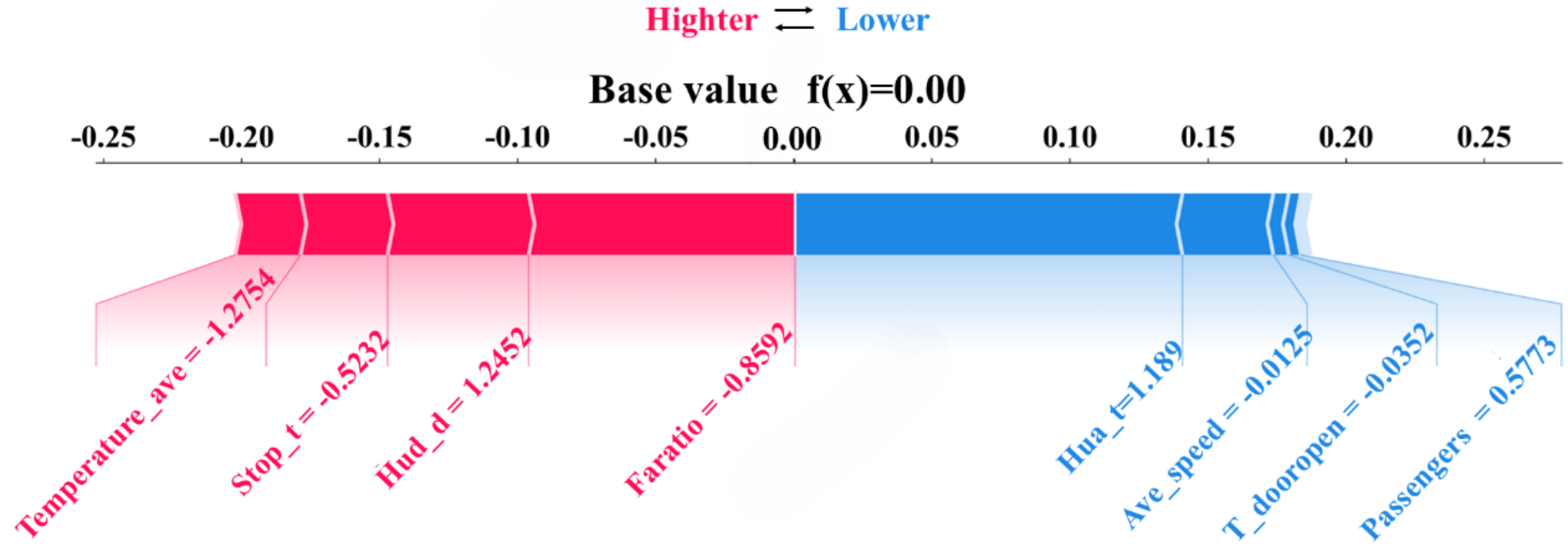

Fig. 11. SHAP force plot.

The cumulative widths of the two camps are roughly balanced, reflecting the complex trade-offs in this operating condition: long coasting creates favorable conditions for saving energy, but its benefit is almost completely offset by Faratio and the low-temperature penalty. Although the individual-level result differs slightly from the global pattern due to local road conditions, the overall direction is consistent with the global regularities. The agreement of the two visualizations on the core attributions confirms the internal stability of the model's decision logic.

### 5.6.5 Nonlinear responses and feature interactions

The SHAP dependence plots (Fig. 12(a)-(i)) map single-feature values to their SHAP contributions, revealing three groups of nonlinear response patterns and two classes of important interactions.

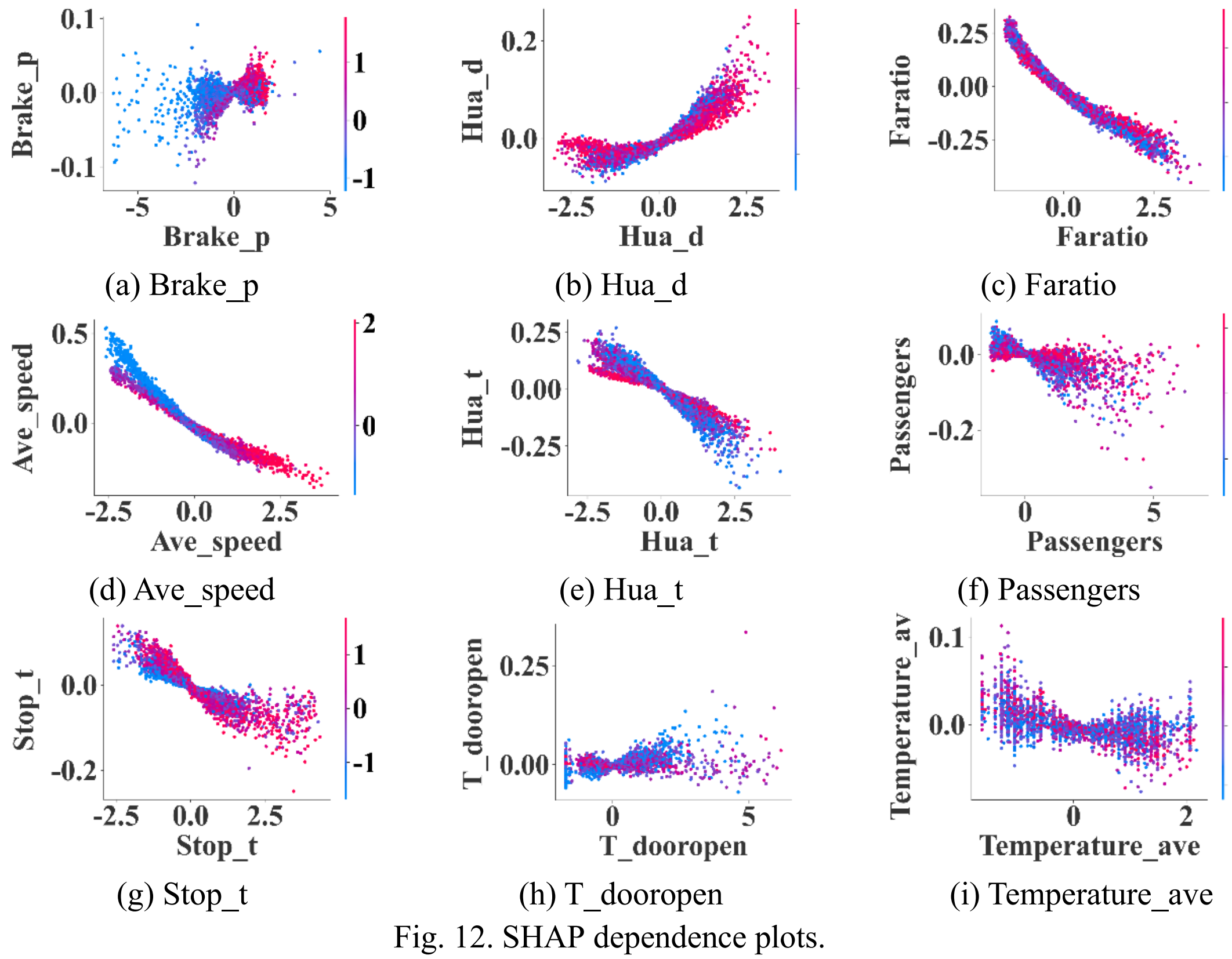

(a) Brake_p (b) Hua_d (c) Faratio
(d) Ave_speed (e) Hua_t (f) Passengers
(g) Stop_t (h) T_dooropen (i) Temperature_ave

Fig. 12. SHAP dependence plots.

For nonlinear responses, the regenerative braking ratio, average speed, and coasting time all correlate negatively with energy use, but the response curves differ in shape. Taking Ave_speed as an example, its SHAP value exhibits a steep positive penalty in the low-speed region (standardized value $< -1$), then drops quickly and flattens in the mid-to-high-speed region, a typical asymmetric threshold effect. This means that lifting bus speed from congested levels to moderate levels brings much larger energy savings than further increases from moderate to high speed. The temperature dependence plot shows a near-linear negative trend, confirming that temperature continuously affects energy use through HVAC load.

For interactions, the color mapping in the dependence plots (the side color bar) reveals important couplings. The partial dependence plot for passenger load reveals clear stratification of scatter points by average speed, suggesting a speed-dependent modulation of the load–energy relationship. Specifically, the energy penalty imposed by passenger weight is amplified under low-speed congestion, whereas its effect is attenuated as vehicle speed increases. Similarly, the joint distribution of average speed and braking time confirms the composite effect of driving behavior: frequent braking at low speeds intensifies energy loss because the recovery potential is limited, while at higher speeds it is partly offset by the larger kinetic-energy-recovery window.

### 5.6.6 Feature interaction network

The feature-interaction network in Fig. 13 quantifies the depth and breadth of coupling between variables from a graph-theoretic perspective. Node color intensity (green gradient) reflects the independent global importance of each feature; edge thickness and color (blue gradient) quantify the interaction strength.

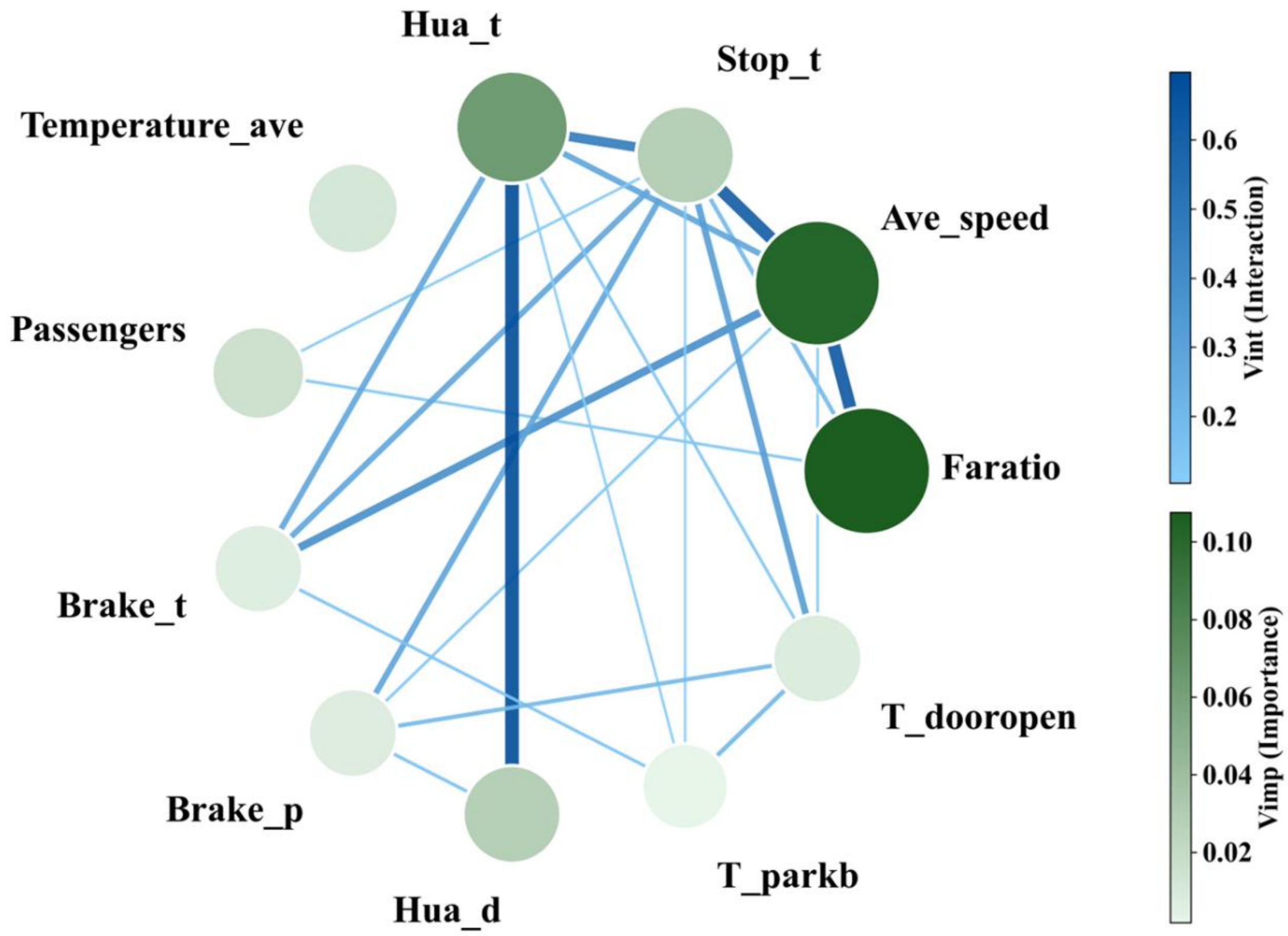

Fig. 13. SHAP-based global interaction network of multi-dimensional operating features.

Three key observations can be drawn from the topology. (1) Coasting time (Hua_t) and coasting distance (Hua_d) show the strongest interaction (0.7), directly mirroring the kinematic identity speed = distance / time, which confirms that the model captures this physical coupling. (2) Dwell time (Stop_t) and average speed (Ave_speed) are strongly coupled, reflecting the typical pattern of low-speed operation combined with frequent starts and stops under congestion. (3) Average speed serves as the dominant hub variable in the feature network, directly associated with regenerative braking ratio, dwell time, and coasting behavior. This centrality is mechanistically grounded in the stop-and-go cycles characteristic of urban bus operations: the operating speed range governs kinetic energy recovery efficiency, while low-speed congestion intensifies start–stop frequency and compresses the coasting window. Through these coupled mechanisms, average speed exerts a broad influence on energy consumption along multiple interacting pathways. Ambient temperature, passenger load, and brake pressure, by contrast, exhibit sparse network connectivity, reflecting relatively independent and less interactive contributions to energy use.

Overall, the interaction network confirms that DTB energy use is genuinely nonlinear: it is not only determined by the individual contributions of each feature but also shaped by the kinematic-state coupling network formed by speed, braking, and coasting.

## 5.7 From correlation to causation

SHAP quantifies marginal contributions but cannot distinguish statistical correlation from causation. This section uses a hierarchical causal-inference framework to decompose energy drivers from both the structural-discovery and effect-estimation perspectives.

### 5.7.1 Correlation analysis of operating and environmental features

Fig. 14 shows the joint evolution of four core operational indicators — monthly mission count, cumulative distance, average passenger load, and energy per km — revealing both periodic rhythms and structural breaks across the operating timeline.

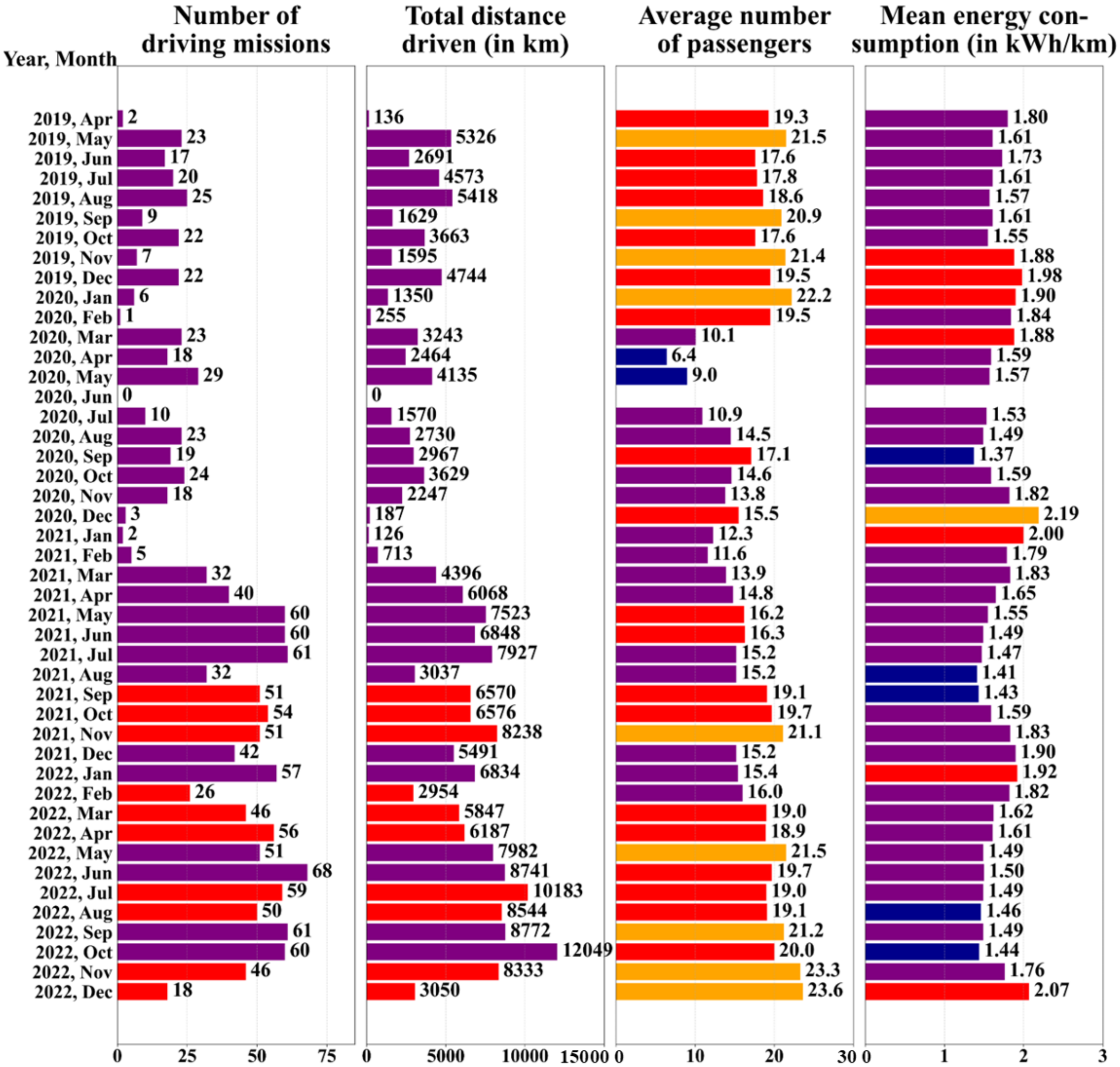

Fig. 14. Range visualization of the ZTBus dataset.

In terms of operating intensity, mission count and total distance change in lock-step. They remain at relatively low levels from 2019 through early 2020, rise sharply from mid-2021, and reach a peak around mid-2022 (some months exceed 60 missions and distances approach 15,000 km), which reflects structural growth driven by network coverage extension and tighter scheduling. The operating trough in the first half of 2020 (with both missions and distance shrinking) aligns with the global decline of public transit during the Coronavirus Disease 2019 (COVID-19) disruption.

The passenger data show two anomalous patterns. (1) From June to July 2020 the average passenger load drops below 10 (dark blue annotation), compared with the typical 15-20, a steep decline directly caused by pandemic-control measures. (2) The orange-annotated high-passenger months cluster around May 2019, November 2021, and May 2022, where the average load exceeds 20, reflecting peaks from holiday surges or moments of tight capacity. As Fig. 15 shows, monthly total ridership fell from about 150,000 to near zero during the pandemic, recovered slowly from early 2021, and gradually returned to historical highs in 2022. This wide, non-stationary distribution of passenger loads introduces diverse and realistic operating conditions into the energy modeling.

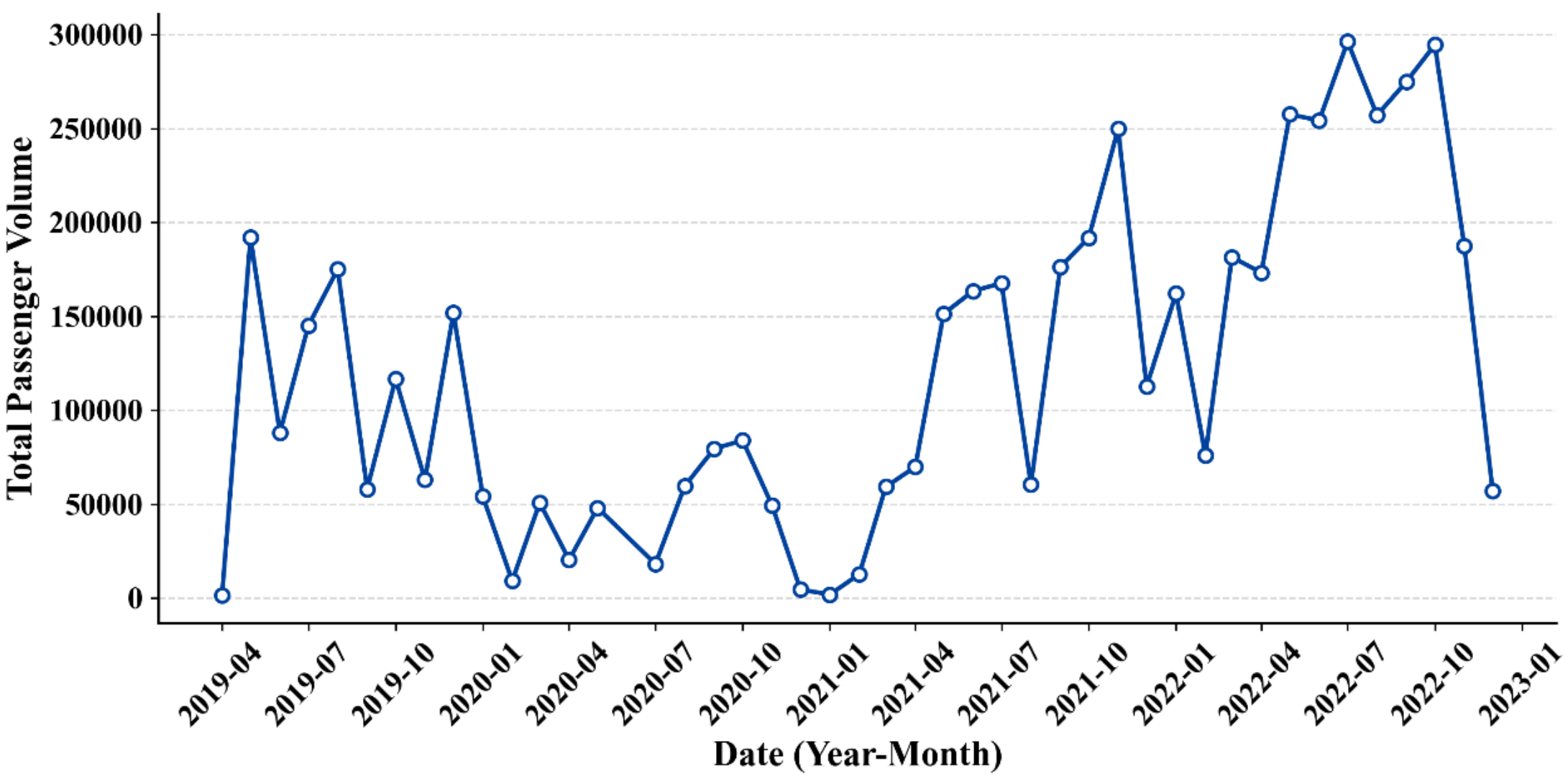

Fig. 15. Monthly passenger volume.

Energy per km follows a clear high-in-winter, low-in-summer rhythm. In winter, low temperatures reduce battery efficiency and increase heating load; some months (for example, December 2020 at 2.19 kWh/km and December 2022 at 2.07 kWh/km) are flagged in orange or red, forming the upper tail of the distribution. In summer and autumn, with mild temperatures and stable auxiliary loads, energy use converges (for example, September 2020 at 1.37 kWh/km and September 2021 at 1.43 kWh/km, flagged in blue). This thermodynamic seasonality, combined with the variation in operating intensity and passenger load, creates the multi-dimensional perturbation structure for energy modeling and provides rich variability for the subsequent causal analysis.

### 5.7.2 Structural path identification with LiNGAM

To directly identify causal directions from observational data, DirectLiNGAM is used to build the DAG over features and energy, shown in Fig. 16. Based on the domain prior of energy conservation, the energy variable is fixed as the terminal node (no outgoing edges), receiving only upstream causal inputs, which rules out the physical contradiction of energy use reversely driving driving behavior.

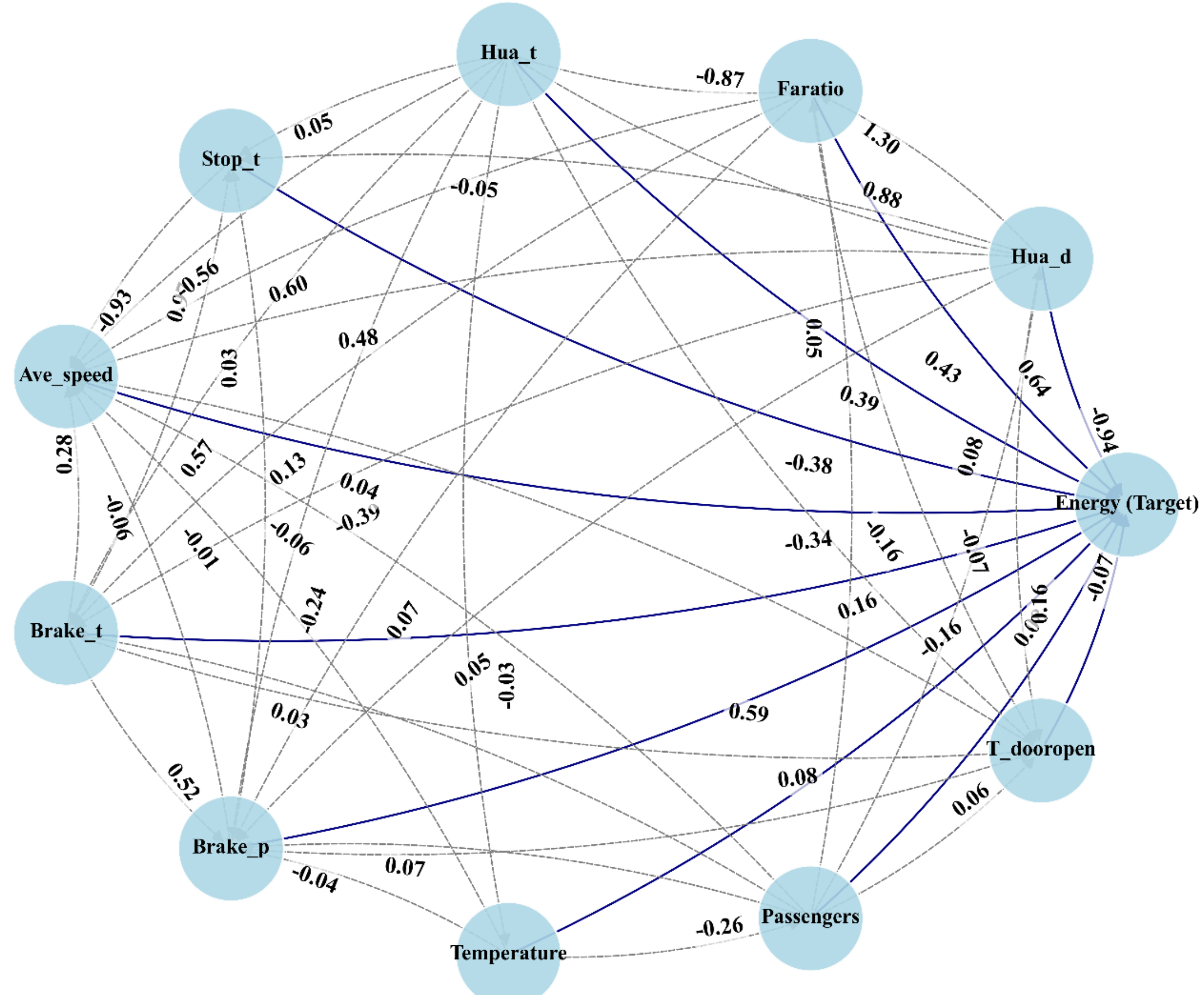

Fig. 16. Directed acyclic graph.

The DAG identifies the regenerative braking ratio (Faratio) and average speed (Ave_speed) as direct causal parents of energy, shown as directed edges from these two variables to energy. From a causal standpoint, the effect of these two variables on energy use is not a statistical co-variation but a directional upstream driver, which cross-validates the SHAP global importance ranking in Section 5.6 from a methodological perspective.

### 5.7.3 Nonlinear intervention-effect estimation with T-Learner

Once the DAG establishes causal directions, the T-Learner meta-learning framework is used to quantify the net causal effect of each driver. Continuous variables are discretized into binary treatment variables using the 75th percentile as the threshold, and the conditional expectations of the treatment and control groups are fitted separately to estimate the ATE. As Fig. 17 shows, the ATEs of the factors differ markedly in direction and magnitude, and can be grouped into energy-saving and energy-increasing drivers.

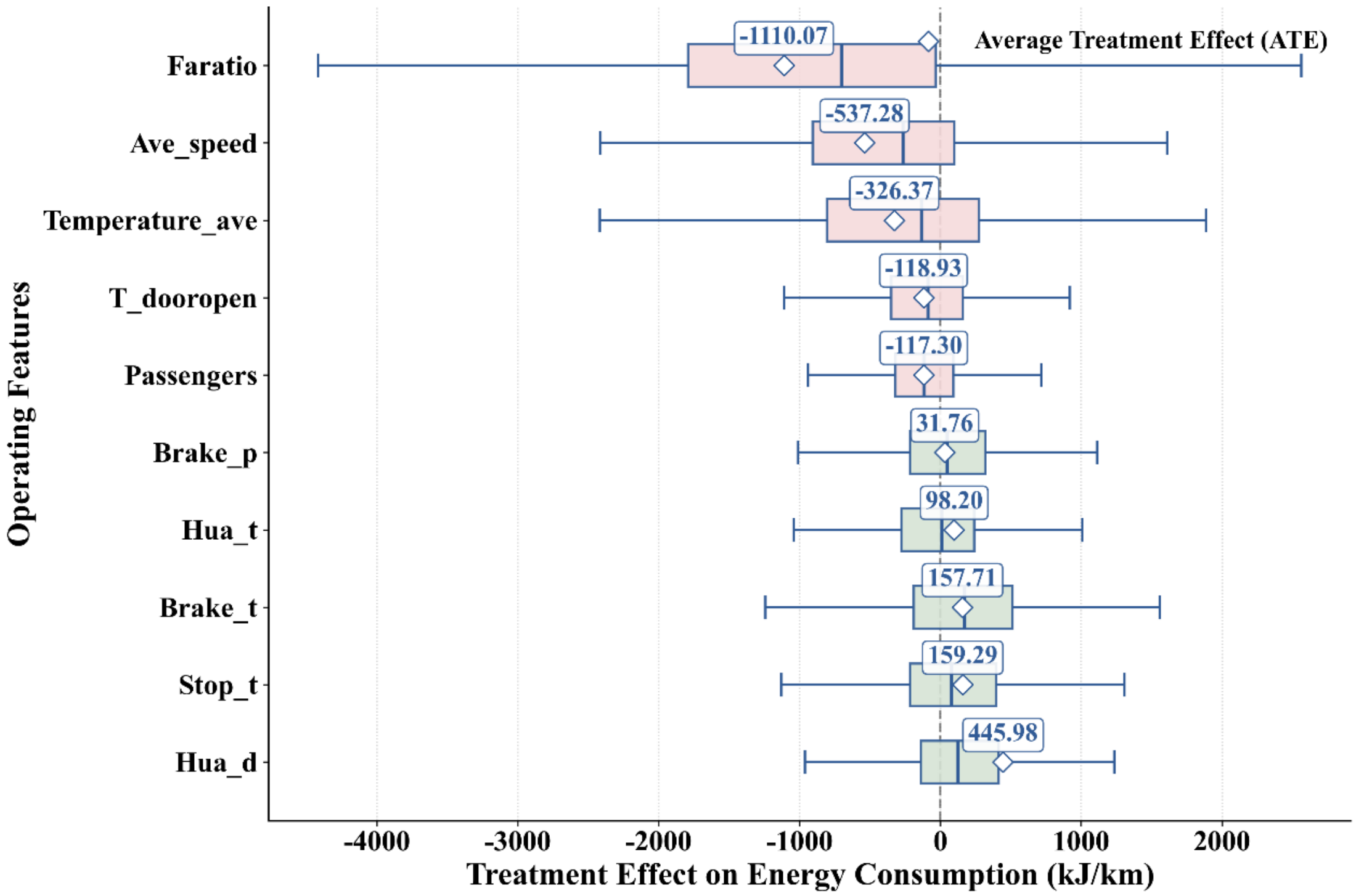

Fig. 17. T-Learner effect quantification.

For the energy-saving drivers, the regenerative braking ratio has the strongest negative causal effect (ATE = −1110.07 kJ/km), indicating that, after controlling for confounders, raising the kinetic-energy recovery ratio is the single most effective intervention for reducing per-kilometer energy use. Average speed also has a strong negative effect (ATE = −537.28 kJ/km), quantifying the energy savings available when moving from low-speed congestion to smoother flow. Combined with the DAG, these two variables form the main intervention anchors in the energy causal chain. The negative effect of passenger load (ATE = −117.30 kJ/km) may seem to contradict the intuitive view that more mass means more energy use, but it is consistent with the energy-recovery dynamics of electric drives: extra load amplifies total inertia, which means more kinetic energy to recover during braking and therefore more electrical energy for the regenerative system. As recent electric-vehicle studies have shown [44], this inertia-driven recovery gain is, under certain conditions, enough to offset the additional traction energy required during acceleration.

For the energy-increasing drivers, coasting distance shows the strongest positive effect (ATE = +445.98 kJ/km). Physically, kinetic energy during coasting is dissipated irreversibly as heat through friction and air drag, cannot be recovered, and must be replaced by new electrical input during the next acceleration, causing a double energy loss. The positive effects of dwell time (ATE = +159.29 kJ/km) and brake time (ATE = +157.71 kJ/km) reflect the same dissipation process from another angle: frequent mechanical braking and operational stagnation both directly consume kinetic energy and indirectly raise overall energy use by prolonging the operation of auxiliary systems in non-traction states.

In summary, the causal-inference results agree with the SHAP analysis on the core findings, but they move the explanation from correlation to mechanism. SHAP reveals which features are important, DirectLiNGAM confirms whether their influence is causal and in which direction, and the T-Learner quantifies the net effect of intervention. Together, they form a progressive chain of explanation and provide a statistically rigorous and physically interpretable basis for operational decisions.

## 5.8 Differentiated contributions of meteorological features

SHAP and causal inference both show that temperature is the most influential meteorological feature, while relative humidity, snowfall, and precipitation contribute relatively little directly. This is driven mainly by the following reasons.

(1) Temperature affects vehicle operation through multiple direct channels. It influences battery charge/discharge efficiency, tire-road friction, HVAC load, and the thermal efficiency of auxiliary systems, and therefore stands out in SHAP importance. Humidity and wind speed, by contrast, usually affect the vehicle through nonlinear interactions or context-specific indirect paths, and their main effects are absorbed by strongly correlated variables such as temperature in a global model.

(2) The variability of the data also matters. The dataset spans a long period, so temperature varies across a wide range that covers all four seasons, while precipitation and wind speed occur less frequently or are concentrated in narrow intervals, so the model naturally picks up the more stable temperature signal. The causal analysis further confirms that temperature has relatively stable direct causal paths to key vehicle indicators.

(3) There is also a nonlinear absorption effect: other meteorological features may be diluted by correlated variables such as time and average speed, making their contribution hard to isolate. Overall, the dominance of temperature in the multi-source dataset reflects both its real physical impact on vehicle operation and the possibility that other meteorological features could be better captured through interaction terms or scenario-specific modeling in future work.

## 5.9 Implications

Based on the experimental results, four sets of implications are drawn.

### (1) Vehicle technology and energy management.

The T-Learner analysis confirms that the regenerative braking ratio is the single most effective intervention for reducing DTB energy use (ATE = −1110.07 kJ/km), about twice as strong as the effect of average speed (−537.28 kJ/km). The SHAP interaction network further shows that the regenerative braking ratio and average speed are tightly coupled through a hub-like structure, with the largest joint energy-saving effect in the mid-to-high-speed range. At the vehicle level, improving the recovery efficiency of the regenerative braking system should be a priority, especially in the 30-50 km/h braking range that corresponds to the most frequent deceleration scenarios in urban transit and offers the largest recovery potential. In addition, both the SHAP beeswarm plot and the correlation analysis show that cold conditions (winter months can reach 2.19 kWh/km) continuously increase energy use through HVAC heating load. Operators are therefore advised to pre-condition vehicles when the daily mean temperature falls below 5 °C — for example, by using night-time low-price electricity or dwell-time charging windows to pre-heat the cabin — and to prioritize heat-pump-based HVAC systems in new vehicle procurement, so as to cut the additional energy draw from electric heaters on the traction battery and markedly reduce the auxiliary share of energy use in the heating season.

### (2) Traffic operation and driving-behavior optimization.

The strong negative causal effect of average speed (ATE = −537.28 kJ/km), combined with the asymmetric threshold pattern in the SHAP dependence plot (a steep positive penalty at low speeds and small marginal gains at higher speeds), indicates that reducing low-speed operation under congestion is a key energy-saving lever. Building on the mapping between average speed and energy consumption, this study proposes segment-specific operational recommendations. For suburban or BRT corridors, the operating speed should be maintained within 40–60 km/h. This interval coincides with the plateau in the SHAP negative-contribution region, where energy consumption reaches its minimum. For inner-city corridors comparable to the ZTBus service area in Zurich, transit signal priority, dedicated bus lanes, and dynamic dispatching are recommended. These measures can raise the segment-level average speed to at least 20–25 km/h. Sustained operation below 15 km/h should be avoided, as the associated energy penalty is substantial. Operators can set bus-priority windows at intersections and congested segments to reduce low-speed dwelling.

Meanwhile, coasting distance shows the strongest positive effect among all energy-increasing factors (ATE = +445.98 kJ/km), and the positive effects of dwell time (+159.29 kJ/km) and brake time (+157.71 kJ/km) are also non-trivial. Together, they indicate that excessive traction-free coasting and frequent mechanical braking are not good energy-saving strategies.Eco-driver training programs should incorporate the following operational guidelines. Coasting distance per event should be kept within 50–150 m; runs exceeding 200 m are discouraged, as prolonged coasting leads to excessive dissipation of kinetic energy. Deceleration should begin approximately 200 m before each scheduled stop, with regenerative braking prioritized to maximize energy recovery. The proportion of mechanical braking events should be kept below 30% of all braking instances. Additionally, single-stop dwell time—encompassing both door-opening and door-closing phases—should be actively managed to curtail unnecessary auxiliary-system energy consumption during stationary periods.

### (3) Service level and capacity allocation.

The negative causal effect of passenger load (ATE = −117.30 kJ/km) carries a policy-relevant finding. As discussed in Section 5.7.3, the underlying mechanism is that additional load amplifies vehicle inertia, which feeds more recoverable energy back into the regenerative braking system and, under suitable conditions, outweighs the extra traction energy required during acceleration. The SHAP dependence plot further shows that this negative effect is stronger in the higher-speed range and weakens under low-speed congestion, where the weight penalty becomes dominant. Clear boundary conditions therefore apply: when the section-average speed exceeds 25 km/h and the regenerative braking system is functioning normally, moderately increasing passenger density does not significantly raise per-kilometer energy use, and in fact favors energy efficiency. For lines meeting these conditions, moderately extending headways during peak hours—for instance, from 3 to 5 minutes—is justified, as it increases per-vehicle passenger load and improves energy efficiency without fundamentally compromising service levels. However, real-time ridership monitoring should be implemented in parallel to ensure that waiting times remain within acceptable thresholds. This safeguard is essential: an unmanaged decline in service quality risks diverting passengers toward private vehicles, thereby undermining the broader sustainability objectives of the transit network. On sections with average speeds below 20 km/h, the extra traction energy from added weight outweighs the regenerative gains, and operators should instead maintain or increase service frequency to spread the load.

### (4) **Policy and planning.**

The three-layer progressive framework — SHAP attribution, DirectLiNGAM causal discovery, and T-Learner effect estimation — provides a systematic quantitative toolkit for network planning and operation. The two intervention anchors identified by the causal analysis (regenerative braking ratio and average speed) can be incorporated into the operational KPI system as key indicators for line-level energy assessment, with suggested quantitative thresholds: line-level regenerative braking ratio no less than 35%, section average speed no less than 20 km/h, and energy use no higher than 2.0 kWh/km.

The causal topology produced by DirectLiNGAM can further help identify energy-sensitive segments on new routes (for example, the recovery potential of long downhill sections and the energy-increasing risk of congested bottlenecks), supporting the layout of catenary networks, charging-facility siting, and supply-coverage design on a causal basis. The core gain from the Time2Vec module confirmed by the ablation study (MAPE reduced from 7.68% to 6.68%) provides a reference for model selection in real-time dispatch systems: models with time-awareness should be prioritized in dispatch systems so that daily and weekly periodic information can be fully exploited to improve scheduling accuracy.

## 6. Conclusion

To address the challenge of accurate and interpretable energy prediction for DTBs, this study builds a unified framework based on Time-Aware TabM and DirectLiNGAM-T-Learner. The main conclusions are as follows.

(1) Combining Time2Vec periodic encoding with the flattened GEMM batch ensemble and tuning hyperparameters by TPE-based Bayesian optimization, Time-Aware TabM reaches a MAPE of 6.52% and an $R^2$ of 0.982, improving MAPE by 1.07 % over LSTM and 1.37 % over XGBoost, and outperforming all ten baselines. The ablation study shows that Time2Vec is the main driver of accuracy (MAPE from 7.68% to 6.68%), while flattened GEMM provides a complementary gain on top of it (down to 6.52%). The MAPE standard deviation (0.07) is also the lowest among all variants, indicating both high accuracy and strong stability.

(2) SHAP multi-layer attribution reveals the dominant role of vehicle operating state features (57.05% of total contribution), identifies the asymmetric threshold effect of average speed in the low-speed region, and uncovers the kinematic coupling network formed by speed, braking, and coasting. Individual attribution confirms that the global patterns also hold at the sample level.

(3) DirectLiNGAM confirms that the regenerative braking ratio and average speed are the direct causal parents of energy, and the T-Learner quantifies the net effects: the regenerative braking ratio (ATE = −1110.07 kJ/km) and average speed (−537.28 kJ/km) are the two main energy-saving drivers, while coasting distance (+445.98 kJ/km) and dwell time (+159.29 kJ/km) are the main energy-increasing factors. The negative effect of passenger load (−117.30 kJ/km) reveals a physical buffering mechanism: greater vehicle mass stores more kinetic energy, a larger share of which can be recovered through regenerative braking.

The study has several limitations that can be addressed in future work. (1) Inference time is about 2.8 times the baseline, so deployment in real-time dispatch still requires model distillation or quantization to reduce latency. (2) The model is validated on only two DTBs in Zurich; future work should extend to more vehicle types and urban environments to assess cross-domain generalization. (3) The linearity assumption of DirectLiNGAM and the binarization used by the T-Learner may discard nonlinear causal information; non-parametric causal-discovery methods or continuous-treatment-effect estimation frameworks could be introduced in future studies. (4) The causal graph is estimated from pooled operational records and represents a single static structure across all routes and seasons. Causal relationships — particularly those involving driving behavior and meteorological conditions — may shift substantially across seasonal contexts or route-specific operating environments. Time-varying or context-stratified causal discovery frameworks would better reflect the structural heterogeneity inherent in real-world transit operations.

# Appendix A

Table A1. Comparison between representative studies and the present work.

| Author | Focus | Method | Interpretability | Findings | Difference |
| --- | --- | --- | --- | --- | --- |
| Fiori et al. (2021) [11] | Microscopic energy modeling of e-buses | Physics-based | High (white-box) | Decouples traction and auxiliary power | Relies on parameter calibration; no weather fusion or temporal encoding |
| Pan et al. (2024) [13] | Online cumulative energy prediction | Longitudinal dynamics + FFRLS | Medium | Stop-by-stop online prediction | No time-awareness or batch ensembling |
| Ferdaus et al. (2026) | Foundation models for clean energy | Review | Low | Zero-shot learning potential | Review only; no DTB modeling |
| [15] | | | | | |
| Lajunen (2014) [14] | Hybrid/BEV bus cost-benefit | Vehicle-dynamics simulation | High | Economic feasibility assessment | Economic focus; no inter-stop prediction |
| Hjelkrem et al. (2022) [12] | Urban bus energy prediction | Longitudinal dynamics + HVAC | High (white-box) | Quantifies HVAC impact | Cannot absorb random traffic disturbances |
| Kim et al. (2025) [16] | Auxiliary energy under extreme climate | Data-driven regression | Low | Quantifies abnormal thermal loads | No driving behavior or temporal features |
| Li et al. (2021) [17] | Trip-level energy decomposition | KNN-RF hybrid | Low | Quantifies temperature/road/traffic impact | No temporal encoding or causal inference |
| Zhang et al. (2020) [18] | EV taxi energy prediction | XGBoost + segment split | Low | Three-level segment framework | Taxi-oriented; no weather fusion or causal analysis |
| Chen et al. (2021) [21] | Instantaneous energy prediction | LSTM + data reconstruction | Low | Real-time prediction | No tabular encoding or causal inference |
| Pan et al. (2025) [22] | Online energy estimation | Incremental learning | Low | Dynamic parameter update | No interpretability or Bayesian optimization |
| Zhu et al. (2025) [23] | Carbon-efficient bus causal analysis | Causal AI framework | High (causal) | Overcomes correlation bias | No deep prediction or multi-source fusion |
| He et al. (2026) [24] | Inter-stop causal effects | Transformer + DML | High (causal) | Estimates causal effects | Different backbone; no TabM or Time2Vec |
| Song et al. (2026) [26] | Interpretable power-load prediction | MSN-KAN | Medium | Enhanced interpretability | Non-transport domain |
| Eskandari et al. (2024) [29] | National energy forecasting | Feature selection + interpretable ML | Medium | Multi-model comparison | Macro scale; no inter-stop prediction |
| Zhang et al. (2025) | Charging-station load | Spatio-temporal causal + game | High (causal) | Disentangles heterogeneous | Charging-station oriented, not bus |

| [30] | balancing | theory | | correlations | energy |
|---|---|---|---|---|---|
| This study | DTB inter-stop energy prediction and causal mechanism | Time-Aware TabM + DirectLiNGAM + T-Learner | High (SHAP + causal) | MAPE 6.52%, ATE quantification | Unified time-aware prediction, interpretability, and causal inference |

Table A2. Fields of the fused inter-stop dataset.

| NO. | Category | Variable | Shorthand | Unit | NO. | Category | Variable | Shorthand | Unit |
|---|---|---|---|---|---|---|---|---|---|
| 1 | Time characteristics | Year | Year | - | 33 | | Traction downhill distance | Traction_d_down | m |
| 2 | | Month | Month | - | 34 | Characteristics of the gliding stage | Coasting energy consumption (theoretical) | Hua_E | kJ |
| 3 | | Date | Day | - | 35 | | Glide time | Hua_t | s |
| 4 | | Hours | Hour | - | 36 | | Glide maximum speed | Hua_maxspeed | m/s |
| 5 | | Sampling interval | Time_diff_s | s | 37 | | Glide average speed | Hua_avespeed | m/s |
| 6 | Operation characteristics | Total trip time | Travel_t | s | 38 | | Coasting distance | Hua_d | m |
| 7 | | Passengers | Passengers | - | 39 | | Standard deviation of glide speed | Hua_speedstd | - |
| 8 | | Catenary supply duration | T_grida | s | 40 | | Glide uphill distance | Hua_d_up | m |
| 9 | | Opening time | T_dooropen | s | 41 | | Distance glided downhill | Hua_d_down | m |
| 10 | | Duration of stops | Stop_t | s | 42 | Braking phase characteristics | Regenerative braking energy | Brake_E | kJ |
| 11 | | Vehicle number | Busnono | - | 43 | | Braking duration | Brake_t | s |
| 12 | Line characteristics | Average speed | Ave_speed | m/s | 44 | | Braking distance | Brake_d | m |
| 13 | | Maximum speed | Max_speed | m/s | 45 | | Braking average deceleration | Brake_aveacc | m/s² |
| 14 | | Speed standard deviation | Std_speed | - | 46 | | Braking uphill distance | Brake_d_up | m |
| 15 | | Travel distance | Travel_d | m | 47 | | Braking downhill distance | Brake_d_down | m |
| 16 | | Uphill distance | D_up | m | 48 | | Pedal pressure | Brake_p | N |
| 17 | | Downhill distance | D_down | m | 49 | | Parking brake duration | T_parkb | s |
| 18 | | Cumulative uphill | Altitude_up | m | 50 | | Regenerative braking ratio | Faratio | % |
| 19 | | Cumulative downhill | Altitude_down | m | 51 | Weather characteristics | Meteorological temperature | Temp | °C |
| 20 | Traction stage characteristics | Traction force | Traction_F | N | 52 | | Precipitation | Precip | mm |
| 21 | | Traction power | Traction_Fv | w | 53 | | Relative humidity | Humidity | % |
| 22 | | Traction time | Traction_t | s | 54 | | Snowfall | Snow | cm |
| 23 | | Traction energy consumption | Traction_E | kJ | 55 | | Snow thickness | Snowdepth | m |
| 24 | | Traction maximum speed | Traction_maxspeed | m/s | 56 | | Wind speed | Windspeed | km/h |
| 25 | | Traction | Traction_avespe | m/s | 57 | | Wind direction | Winddir | ° |

| | | average speed | ed | | | | | | |
|---|---|---|---|---|---|---|---|---|---|
| 26 | | Traction distance | Traction_d | m | 58 | | Sea level pressure | Sealevelpressure | kPa |
| 27 | | Traction average acceleration | Traction_aveacc | m/s² | 59 | | Average temperature in the car | Temperature_ave | °C |
| 28 | | Traction maximum acceleration | Traction_max_acc | m/s² | 60 | | Energy consumption per mileage | Energy_per_km | kJ/km |
| 29 | | Traction minimum acceleration | Traction_min_acc | m/s² | 61 | Energy consumption characteristics | Maximum energy consumption in the range | Max_e | kJ |
| 30 | | Standard deviation of traction speed | Traction_speedstd | - | 62 | | Crowding ratio | Fratio | % |
| 31 | | Standard deviation of traction acceleration | Traction_std_acc | - | 63 | | Load ratio | Weightratio | % |
| 32 | | Traction uphill distance | Traction_d_up | m | | | | | |

## Data availability

The datasets generated and analysed during the current study will be made publicly available via a persistent URL upon formal publication of the article, in accordance with funder and institutional requirements for data sharing.

## Credit authorship contribution statement

**Wentao Zeng:** Investigation, Data curation, Validation, Visualization, Writing-Original draft, Methodology. **Zijian Huang:** Supervision, Writing-Reviewing, Methodology. **Yiming Bie:** Supervision, Funding acquisition. **Jiabin Wu:** Conceptualization, Methodology, Funding acquisition, Supervision, Writing-Reviewing and Editing. **Jun Gong:** Investigation.

## Declaration of Competing Interest

The authors declare that they have no known competing financial interests or personal relationships that could have appeared to influence the work reported in this paper.

## Acknowledgments

This research is supported by the Guangdong Province Philosophy and Social Science Planning Youth Project [grant number GD26YGG27]; Natural Science Foundation of Guangdong Province, China [grant number 2026A1515011378]; National Natural Science Foundation of China [grant number 72471102];

## References

[1] Zhang H, Wang S, Hao J, et al. Air pollution and control action in Beijing. Journal of Cleaner Production, 2016, 112: 1519-1527. DOI: 10.1016/j.jclepro.2015.04.092.

[2] Tamayao M-A M, Michalek J J, Hendrickson C, et al. Regional variability and uncertainty of electric vehicle life cycle CO2 emissions across the United States. Environmental Science & Technology,

2015, 49(14): 8844-8855. DOI: 10.1021/acs.est.5b00815.
[3] Hu B, Gao C, Gao W. Integrated optimization of timetabling and berthing for dual-source trolleybus routes along corridors. Computers & Industrial Engineering, 2025, 203: 111066. DOI: 10.1016/j.cie.2025.111066.
[4] IEA. Global EV Outlook 2025. https://www.iea.org/reports/global-ev-outlook-2025.
[5] Wu J, Huang Z, Zhou W. Energy consumption prediction for dual-powered trolley bus based on a SE-TCN hybrid network. Measurement, 2026, 262: 119978. DOI: 10.1016/j.measurement.2025.119978.
[6] Lisicki M, Lubitz W, Taylor G W. Optimal design and operation of Archimedes screw turbines using Bayesian optimization. Applied Energy, 2016, 183: 1404-1417. DOI: 10.1016/j.apenergy.2016.09.084.
[7] Hajihosseinlou M, Maghsoudi A, Ghezelbash R. Development of a hybrid Bayesian-Hyperband optimization procedure: GeoAI-driven hyperparameter tuning of AdaBoost for enhancing mineral prospectivity mapping. Environmental Modelling and Software, 2026, 198: 106883. DOI: 10.1016/j.envsoft.2026.106883.
[8] Meng F, Qiao S, Chen D. A data-driven framework for precise geometric measurement of tunnel structures using 3D point clouds and Bayesian optimization. Measurement, 2026: 120931. DOI: 10.1016/j.measurement.2026.120931.
[9] Fu Z, Yang X, Ma Y, et al. Integrating explainable AI and causal inference to unveil regional air quality drivers in China. Journal of Environmental Management, 2025, 390: 126270. DOI: 10.1016/j.jenvman.2025.126270.
[10] Bareinboim E, Pearl J. Causal inference and the data-fusion problem. Proceedings of the National Academy of Sciences, 2016, 113(27): 7345-7352. DOI: 10.1073/pnas.1510507113.
[11] Fiori C, Montanino M, Nielsen S, et al. Microscopic energy consumption modelling of electric buses: model development, calibration, and validation. Transportation Research Part D: Transport and Environment, 2021, 98: 102978. DOI: 10.1016/j.trd.2021.102978.
[12] Hjelkrem O A, Lervåg K Y, Babri S, et al. A battery electric bus energy consumption model for strategic purposes: validation of a proposed model structure with data from bus fleets in China and Norway. Transportation Research Part D: Transport and Environment, 2021, 94: 102804. DOI: 10.1016/j.trd.2021.102804.
[13] Pan Y, Fang W, Ge Z, et al. A hybrid on-line approach for predicting the energy consumption of electric buses based on vehicle dynamics and system identification. Energy, 2024, 290: 130205. DOI: 10.1016/j.energy.2023.130205.
[14] Lajunen A. Energy consumption and cost-benefit analysis of hybrid and electric city buses. Transportation Research Part C: Emerging Technologies, 2014, 38: 1-15. DOI: 10.1016/j.trc.2013.10.008.
[15] Ferdaus M M, Dam T, Sarkar M R, et al. Foundation models for clean energy forecasting: a comprehensive review. Renewable and Sustainable Energy Reviews, 2026, 226: 116452. DOI: 10.1016/j.rser.2025.116452.
[16] Kim D, Yun J, Jang K, et al. Auxiliary energy consumption of electric vehicles: modeling and prediction using real-world vehicle data. Applied Energy, 2025, 401: 126766. DOI: 10.1016/j.apenergy.2025.126766.
[17] Li P, Zhang Y, Zhang Y, et al. Prediction of electric bus energy consumption with stochastic speed

profile generation modelling and data-driven method based on real-world big data. Applied Energy, 2021, 298: 117204. DOI: 10.1016/j.apenergy.2021.117204.
[18] Zhang J, Wang Z, Liu P, et al. Energy consumption analysis and prediction of electric vehicles based on real-world driving data. Applied Energy, 2020, 275: 115408. DOI: 10.1016/j.apenergy.2020.115408.
[19] Zhang Z, Wang R, Liu P, et al. Research on energy consumption law and charging strategies design of electric buses. Energy, 2025, 322: 135327. DOI: 10.1016/j.energy.2025.135327.
[20] Liu Z, Song Z. Robust planning of dynamic wireless charging infrastructure for battery electric buses. Transportation Research Part C: Emerging Technologies, 2017, 83: 77-103. DOI: 10.1016/j.trc.2017.07.013.
[21] Chen Y, Zhang Y, Sun R. Data-driven estimation of energy consumption for electric bus under real-world driving conditions. Transportation Research Part D: Transport and Environment, 2021, 98: 102969. DOI: 10.1016/j.trd.2021.102969.
[22] Pan Y, Fang W, Yan H, et al. An online energy consumption estimation method for different types of battery electric buses based on incremental learning. Energy, 2025, 336: 138515. DOI: 10.1016/j.energy.2025.138515.
[23] Zhu B, Shu S, Roncoli C, et al. Causal AI for carbon-efficient bus systems: unveiling the impacts of policy interventions. Transportation Research Part D: Transport and Environment, 2025, 149: 105034. DOI: 10.1016/j.trd.2025.105034.
[24] He Y, Huang Z, Zhan J, et al. A transformer-enhanced causal analysis framework for inter-stop energy consumption of electric bus. Transportation Research Part D: Transport and Environment, 2026, 153: 105221. DOI: 10.1016/j.trd.2026.105221.
[25] Llopis-Castelló D, Arce-Blanco C, Valdecantos-Álvarez J C. Before-after study of pavement roughness impact on vehicle fuel consumption and emissions. Transportation Research Part D: Transport and Environment, 2026, 151: 105131. DOI: 10.1016/j.trd.2025.105131.
[26] Song L, Wang B, Liu Y, et al. Enhancing power load forecasting accuracy under high renewable energy penetration with a MSN-KAN framework. Results in Engineering, 2026, 29: 109275. DOI: 10.1016/j.rineng.2026.109275.
[27] Zhang Y, Li H, Ren G. Quantifying the social impacts of the London Night Tube with a double/debiased machine learning based difference-in-differences approach. Transportation Research Part A: Policy and Practice, 2022, 163: 288-303. DOI: 10.1016/j.tra.2022.07.015.
[28] Jia S, Liu X, Zhao L, et al. Interpretable spatiotemporal urban energy forecasting. Energy, 2025, 334: 137503. DOI: 10.1016/j.energy.2025.137503.
[29] Eskandari H, Saadatmand H, Ramzan M, et al. Innovative framework for accurate and transparent forecasting of energy consumption: a fusion of feature selection and interpretable machine learning. Applied Energy, 2024, 366: 123314. DOI: 10.1016/j.apenergy.2024.123314.
[30] Zhang X, Zhao Q, Wang S. A spatiotemporal causal inference based multi-level game framework for load balancing among fast charging stations. Energy, 2025, 336: 138394. DOI: 10.1016/j.energy.2025.138394.
[31] Widmer F, Ritter A, Onder C H. ZTBus: a large dataset of time-resolved city bus driving missions. Scientific Data, 2023, 10(1): 687. DOI: 10.1038/s41597-023-02600-6.

[32] Zeng X, Wu J, Ling M, et al. Dual-powered trolley bus inter-stop travel energy consumption prediction method: temporal feature transformer framework. Transportation Research Record, 2026, 2680(2): 580-597. DOI: 10.1177/03611981251366263.

[33] Gorishniy Y, Kotelnikov A, Babenko A. TabM: advancing tabular deep learning with parameter-efficient ensembling. https://github.com/yandex-research/tabm.

[34] Shimizu S, Inazumi T, Sogawa Y, et al. DirectLiNGAM: a direct method for learning a linear non-Gaussian structural equation model. arXiv, 2011. DOI: 10.48550/arXiv.1101.2489.

[35] Box G E P, Jenkins G M. Time series analysis: forecasting and control. The Statistician, 1978, 27(3/4): 265. DOI: 10.2307/2988198.

[36] Drucker H, Burges C J C, Kaufman L, et al. Support vector regression machines. In: Advances in Neural Information Processing Systems, 1996, 9.

[37] Breiman L. Random forests. Machine Learning, 2001, 45(1): 5-32. DOI: 10.1023/A:1010933404324.

[38] Chen T, Guestrin C. XGBoost: a scalable tree boosting system. In: Proc. 22nd ACM SIGKDD, 2016: 785-794. DOI: 10.1145/2939672.2939785.

[39] Ke G, Meng Q, Finley T, et al. LightGBM: a highly efficient gradient boosting decision tree. In: Advances in Neural Information Processing Systems, 2017, 30.

[40] Pamuła T, Pamuła D. Prediction of electric buses energy consumption from trip parameters using deep learning. Energies, 2022, 15(5): 1747. DOI: 10.3390/en15051747.

[41] LeCun Y, Bottou L, Bengio Y, et al. Gradient-based learning applied to document recognition. Proceedings of the IEEE, 1998, 86(11): 2278-2324. DOI: 10.1109/5.726791.

[42] Hochreiter S, Schmidhuber J. Long short-term memory. Neural Computation, 1997, 9(8): 1735-1780. DOI: 10.1162/neco.1997.9.8.1735.

[43] Hollmann N, Müller S, Purucker L, et al. Accurate predictions on small data with a tabular foundation model. Nature, 2025, 637(8045): 319-326. DOI: 10.1038/s41586-024-08328-6.

[44] Moawad A, Zebiak M, Pierre M S, et al. Insights into adaptive cruise control and energy efficiency in electric vehicles. npj Sustainable Mobility and Transport, 2025, 2(1): 49. DOI: 10.1038/s44333-025-00066-0.